\documentclass[oneside, 12pt]{amsart}
\usepackage[a4paper, top=3cm, left=2.8cm, right=2.8cm]{geometry}
\usepackage{amscd, amssymb, amsmath, mathrsfs}
\usepackage[english]{babel}
\usepackage{url}
\usepackage{tikz-cd}
\usepackage{booktabs}
\usepackage[pdftex, colorlinks=true,  citecolor=blue, linkcolor=blue, linktocpage=true]{hyperref}

\newcommand\mylabel[1]{\label{#1}\marginpar{\vspace{-1ex}\medskip\medskip\footnotesize \tt #1}}
\renewcommand\mylabel[1]{\label{#1}}
\newcommand{\mydate}{
\number\day\space
\ifcase\month \or January\or February\or March\or April\or May\or June\or July\or August\or September\or October\or November\or December\fi 
\space\number\year}

\DeclareUrlCommand\arXiv{\urlstyle{same}}

\newtheorem{theorem}{Theorem}[section]
\newtheorem{maintheorem}{Theorem}

\newtheorem{lemma}[theorem]{Lemma}
\newtheorem{proposition}[theorem]{Proposition}
\newtheorem{corollary}[theorem]{Corollary}

\theoremstyle{definition}
\newtheorem{definition}[theorem]{Definition}

\newtheorem*{acknowledgement}{Acknowledgement}

\theoremstyle{remark}

\newcommand{\ZZ}{\mathbb{Z}}
\newcommand{\QQ}{\mathbb{Q}}
\newcommand{\RR}{\mathbb{R}}
\newcommand{\CC}{\mathbb{C}}
\newcommand{\FF}{\mathbb{F}}

\newcommand{\PP}{\mathbb{P}}
\renewcommand{\AA}{\mathbb{A}}
\newcommand{\GG}{\mathbb{G}}

\newcommand{\ideala}{\mathfrak{a}}
\newcommand{\idealb}{\mathfrak{b}}

\newcommand{\shF}{\mathscr{F}}

\newcommand{\shI}{\mathscr{I}}

\newcommand{\shM}{\mathscr{M}}

\newcommand{\shL}{\mathscr{L}}

\newcommand{\Aff}{\text{\rm Aff}}
\newcommand{\alg}{\text{\rm alg}}

\newcommand{\Ann}{\operatorname{Ann}}

\newcommand{\Aut}{\operatorname{Aut}}

\newcommand{\Bl}{\operatorname{Bl}}

\newcommand{\cl}{\operatorname{cl}}

\newcommand{\disc}{\operatorname{disc}}

\newcommand{\End}{\operatorname{End}}

\newcommand{\Et}{{\text{\rm Et}}}

\newcommand{\Exc}{\operatorname{Exc}}

\newcommand{\Frac}{\operatorname{Frac}}

\newcommand{\Gal}{\operatorname{Gal}}
\newcommand{\GL}{\operatorname{GL}}

\newcommand{\I}{\text{\rm I}}
\newcommand{\II}{\text{\rm II}}

\newcommand{\IV}{\text{\rm IV}}

\newcommand{\Kernel}{\operatorname{Ker}}

\newcommand{\invlim}{\varprojlim}

\newcommand{\lra}{\longrightarrow}

\newcommand{\maxid}{\mathfrak{m}}

\newcommand{\primid}{\mathfrak{p}}
\renewcommand{\O}{\mathscr{O}}

\newcommand{\op}{\text{\rm op}}

\newcommand{\perf}{{\text{\rm perf}}}

\newcommand{\PGL}{\operatorname{PGL}}

\newcommand{\pr}{\operatorname{pr}}
\newcommand{\Proj}{\operatorname{Proj}}

\newcommand{\quadand}{\quad\text{and}\quad}

\newcommand{\ra}{\rightarrow}

\newcommand{\rank}{\operatorname{rank}}
\newcommand{\red}{{\operatorname{red}}}
\newcommand{\Reg}{\operatorname{Reg}}

\newcommand{\Res}{\operatorname{Res}}

\newcommand{\sep}{{\operatorname{sep}}}

\newcommand{\sh}{{\text{\rm sh}}}

\newcommand{\Sing}{\operatorname{Sing}}
\newcommand{\SL}{\operatorname{SL}}

\newcommand{\Sp}{\operatorname{Sp}}
\newcommand{\Spec}{\operatorname{Spec}}

\newcommand{\uH}{\underline{H}}
\newcommand{\val}{\operatorname{val}}

\newcommand{\bu}{\bar{u}}
\newcommand{\bv}{\bar{v}}

\newcommand{\lieg}{\mathfrak{g}}

\newcommand{\Fix}{\operatorname{Fix}}

\newcommand{\Zero}{\operatorname{Zero}}
\newcommand{\crit}{\text{\rm crit}}
\newcommand{\cyc}{\text{\rm cyc}}
\newcommand{\uRes}{\underline{\operatorname{Res}}}
\newcommand{\univ}{\text{\rm univ}}

\newcommand{\tK}{{\tilde{K}}}
\newcommand{\Kthree}{\text{\rm K3}}
\newcommand{\Enr}{\text{\rm Enr}}
\newcommand{\good}{\text{\rm good}}
\newcommand{\bad}{\text{\rm bad}}

\begin{document}

\title[K3 surfaces over small number fields]
      {K3 surfaces over small number fields  and Kummer constructions in families}

\author[Stefan Schr\"oer]{Stefan Schr\"oer}
\address{Heinrich Heine University D\"usseldorf, Faculty of Mathematics and Natural Sciences, Mathematical Institute, 40204 D\"usseldorf, Germany}
\curraddr{}
\email{schroeer@math.uni-duesseldorf.de}


\dedicatory{Revised version, 5 September 2026}

\begin{abstract}
We construct   K3 surfaces over   number fields that have good reduction everywhere.
These do not exist over the  rational numbers, by results of Abrashkin and Fontaine.
Our surfaces exist  for  three quadratic number fields, and an infinite family of $S_3$-number fields. 
To this end we develop a theory of Kummer constructions in families, based on Romagny's notion of the 
effective models, here  applied to sign involutions. This includes   quotients of non-normal surfaces
by infinitesimal group schemes 
in characteristic two, as developed by Kondo and myself.
By the results of Brieskorn and Artin, the resulting families of normal K3 surfaces admit simultaneous resolutions of singularities,
at least  after suitable base-changes. These resolutions are constructed in two ways:
First, by blowing-up families of one-dimensional centers 
that acquire embedded components. Second, by computing various $\ell$-adic local systems in terms of representation theory,
and invoking Shepherd-Barron's results on the resolution functor, which is representable by   a highly non-separated algebraic space.
\end{abstract}

\maketitle
\tableofcontents

\section*{Introduction}
\mylabel{Introduction}

Besides obvious examples like projective spaces $\PP^n$, toric varieties $\operatorname{Temb}_N(\Delta)$ 
stemming from regular fans, or flag varieties $G/P$
coming from  Chevalley groups, it is apparently  
exceedingly difficult to specify schemes $X$ over the ring  $R=\ZZ$
so that the structure morphism $X\ra\Spec(\ZZ)$ is smooth and proper.
The Minkowski discriminant bound   for  number fields shows that there
are no finite \'etale morphisms $\Spec(\O_F)\ra\Spec(\ZZ)$ of degree at least two 
(\cite{Minkowski 1896}, Section 42 or \cite{Neukirch 1999}, Chapter III, Theorem 2.17).
By deep results of Abrashkin and Fontaine (\cite{Abrashkin 1985} and \cite{Fontaine 1985}), 
there are no families of abelian varieties $A\ra\Spec(\ZZ)$
of relative dimension at least one.
More generally, they independently  established that for any smooth proper $X\ra \Spec(\ZZ)$,
the complex fiber $X_\CC$ satisfies severe restrictions on the Hodge numbers (\cite{Abrashkin 1990} and \cite{Fontaine 1993}).
As a consequence,  there are no families of K3 surfaces over the integers. Building on this, the author
proved that there are no such families of Enriques surfaces either \cite{Schroeer 2023}.

Recall that the defining properties of \emph{K3 surfaces} are $c_1=0$ and $b_2=22$, with   quartic surfaces
in $\PP^3$ as simplest examples.
Given the  extraordinary role of K3 surfaces in algebraic geometry, it is natural to ask \emph{if they exist  in families
over number rings with small degree and discriminant}. In other words, we seek to understand
for which number fields $\QQ\subset F$  the stack of K3 surfaces 
$$
\shM_\Kthree\lra(\Aff/S)
$$ 
starts to have non-empty fiber category $\shM_\Kthree(\O_F)$. This category is  empty for $F=\QQ$, 
and according to Abrashkin  (\cite{Abrashkin 1985}, Theorem in \S 7, Section 6) 
remains empty for $F=\QQ(\sqrt{d})$ with $d\in\{-1,-3,5\}$.
Andr\'e  (\cite{Andre 1996}, Theorem 1.3.1) showed that for fixed  $F$ and  
$R=\O_F[1/a]$,  up to isomorphism there are only finitely many polarized families
$X\ra\Spec(R)$ belonging to  $\shM_\Kthree(R)$,
an instance of the so-called \emph{Shafarevich Problem} \cite{Shafarevich 1963}.
Analogues  of the N\'eron--Ogg--Shafarevich Criterion for good reduction of abelian varieties were
developed  for K3 surfaces by Matsumoto and Liedtke 
\cite{Matsumoto 2015a}, \cite{Matsumoto 2015b}, \cite{Liedtke; Matsumoto 2018}, \cite{Matsumoto 2023} 
and Bragg and Yang \cite{Bragg; Young 2023}, see also  \cite{Chiarellotto; Lazda; Liedtke 2019},  \cite{Chiarellotto; Lazda; Liedtke 2022}. 
Reductions of Kummer surfaces in characteristic two, a case of particular relevance for us, were
studied by Overkamp \cite{Overkamp 2021} and Lazda and Skorobogatov \cite{Lazda; Skorobogatov 2023}.
Our first main  result reveals that the fiber categories for $\shM_\Kthree$ are
non-empty sooner than one perhaps might expect:

\begin{maintheorem}
(See Thm.\ \ref{k3 over quadratic fields} and \ref{k3 over s3 fields})
The fiber category $\shM_\Kthree(\O_F)$ is non-empty if $\QQ\subset F$ is  
\begin{enumerate}
\item 
the $S_3$-number field arising as Galois closure from a cubic number field $K$ whose discriminant takes the form 
$d_K=-3f^2$ for some even factor $f$;
\item 
one of the three  quadratic number fields with discriminant $d_F\in\{28,41,65\}$.
\end{enumerate}
\end{maintheorem}

In turn, the smallest  examples are $F=\QQ(\sqrt[3]{2},e^{2\pi i/3})$ and $F=\QQ(\sqrt{7})$.
The desired families   of K3 surfaces will be constructed from families of normal Kummer surfaces.
In the case of quadratic number fields we also get:

\begin{maintheorem}
(See Thm.\ \ref{enriques over quadratic fields})
The fiber category $\shM_\Enr(\O_F)$ for the  stack of Enriques surfaces is non-empty for the three quadratic
number fields $F$ with discriminant $d_F\in\{28,41,65\}$.
\end{maintheorem}
 
In some sense, the above results are neat by-products. My main motivation is to develop a \emph{theory for    Kummer constructions in families}.
The classical Kummer construction goes back to the nineteenth century:
In modern terms, it attaches to each abelian surface $A$ the quotient $A/\{\pm 1\}$ by the sign involution.
This is a normal surface with rational double points and the minimal
resolution is a K3 surface.
This remains  only partly true in characteristic $p=2$: Shioda \cite{Shioda 1974} and  Katsura \cite{Katsura 1978} observed
that if $A$ is supersingular,  the quotient
by the sign involution   yields rational surfaces rather than K3 surfaces.
To remedy this confusing situation,  the author replaced $A$ by the self-product $C\times C$ of the rational cuspidal curve,
a non-normal genus-one curve, 
and the  action of    $\{\pm 1\}$ by an action of the  infinitesimal group scheme $\alpha_2$, and showed
that $(C\times C)/\alpha_2$ is indeed a K3 surface with rational double points \cite{Schroeer 2007}.
More recently, Kondo and the author \cite{Kondo; Schroeer 2021} considered actions of the multiplicative group scheme
$\mu_2$, and showed that $(C\times C)/\mu_2$ is a K3 surface with only sixteen rational double points
of type $A_1$, as in the classical situation. In both cases, however, there is  one additional $D_4$-singularity.

We now put ourselves into a relative situation, over a Dedekind scheme $S$. Let $E_F$ be an elliptic curve
over the function field $F=\kappa(\eta)$, with resulting Weierstra\ss{} model $E\ra S$.
At the residue fields of characteristic two, we allow  additive reduction, but discard
supersingular good reduction. 
The central idea of this paper is to use Romagny's powerful theory of \emph{effective models of group scheme
actions} \cite{Romagny 2012}, and replace the sign involution by its effective model
$$
G=\overline{\{\pm 1\}_\eta}\subset\Aut_{E/S}.
$$
This is a family of group schemes of order two, acting faithfully on the fibers of the Weierstra\ss{} model,
even at the points of bad reduction. It is then possible to characterize the condition that $G$ acts
freely at  the scheme of non-smoothness $\Sing(E/S)$, the crucial prerequisite for   infinitesimal Kummer constructions.
This leads to the notion of elliptic curves $E_F$ and $E'_F$ that are \emph{admissible for the Kummer construction}
(Definition \ref{admissible for kummer}).
This crucial condition pertains to the primes of bad reduction with residue characteristic two: The reduction 
type must be additive, and the action of the effective model at the singular point is free.
This can be re-formulated in terms of the local Weierstra\ss{} equation, and is equivalent to  the opaque condition
$\val(b_4a_2+b_6)=2\val(2,a_1,a_3)$. Our third main result: 
 
\begin{maintheorem}
(See Thm.\ \ref{quotient is family of k3})
If the elliptic curve  $E_F$   is admissible for the Kummer construction, then the categorical quotient 
$V=(E\times E)/G$ is a family of normal  K3 surfaces over the Dedekind scheme $S$.
\end{maintheorem}

The existence of such categorical quotients easily follows from the Keel--Mori Theorem
\cite{Keel; Mori 1997}, and was further analyzed  by Rydh \cite{Rydh 2013}.
The above also works for pairs $E_F$, $E'_F$ of elliptic curves that are admissible for the Kummer construction.
In the future, we  hope to extend the  theory to general abelian surfaces $A_F$.

Let $g:V\ra S$ be the structure morphism for $V=(E\times E)/G$. Then the singularities in the geometric fibers are at most rational double points.
The  problem now is to construct  a \emph{simultaneous minimal resolution} $r:X\ra V$.
The main insight of Brieskorn (\cite{Brieskorn 1966}, \cite{Brieskorn 1968}, \cite{Brieskorn 1971})    was that this indeed exists, at least 
after making a base-change along some finite surjective $S'\ra S$. Building on this, Artin \cite{Artin 1974} showed that
the functor of resolutions is actually representable by an algebraic space $\Res_{V/S}$, coming with a \emph{universal
simultaneous resolution}
$$
f_\univ:X_\univ\ra\Res_{V/S}
$$
Usually, the base and therefore also the total space  are highly non-separated, in a mind-blowing un-schematic way.
On the other hand, one may argue that the base is rather close to being a scheme, because
the  continuous map $|\Res_{V/S}|\ra |S|$ is bijective. 
In Lemma \ref{purely inseparable extensions} we observe  directly that the residue field extensions 
$\kappa(s)\subset\kappa(s')$ are  purely inseparable. Here $s'\in|\Res_{V/S}|$ signifies the point corresponding to $s\in S$.

Shepherd-Barron \cite{Shepherd-Barron 2021} makes a  detailed analysis of    obstructions in $\ell$-adic cohomology
for the existence of a simultaneous minimal resolution $X\ra V$, and relates them to the Weyl groups $W$
attached to the RDP singularities via dual graphs as  Dynkin diagrams. 
Applying  his abstract results for our concrete  $V=(E\times E')/G$, we are able to pin-point the geometric meaning of the 
obstruction, and thus explain exactly what base-changes are needed.
The relevant part in $\uH^2(X_\eta,\ZZ_\ell)$ is a Kronecker product arising from the K\"unneth Theorem,
which has to match with the part in $\uH^2(X_\sigma,\ZZ_\ell(1))$ stemming from the cycle class map for 
the exceptional divisors of the additional $D_4$-singularities $v_\crit\in V$. Using the representation-theoretic fact
that Kronecker products are almost never permutation representations, we get  sufficient control over the situation:

\begin{maintheorem}
(See Thm.\ \ref{sufficient condition for resolution})
A simultaneous minimal resolution of singularities $r:X\ra V$ exists provided the following conditions hold:
\begin{enumerate}
\item The field $F$ contains a primitive third root of unity.
\item The purely inseparable extension $k\subset \kappa(v_\crit)$ is an equality.
\item The elliptic curves $E\otimes F^\sh$ and $E'\otimes F^\sh$ are isomorphic, and acquire good reduction 
over a quadratic extension $F^\sh\subset L$. 
\end{enumerate}
\end{maintheorem}

An important intermediate step that produces \emph{simultaneous partial resolutions}
is provided by  blowing-ups $Y=\Bl_Z(V)$ with families of one-dimensional schemes, where the closed fiber acquires   embedded points,
a new technique that should be useful in many other contexts. Our main result in this direction
is Theorem \ref{simultaneous partial resolution}.
 
Coming back to the stack of K3 surfaces $\shM_\Kthree$, the problem now is 
to find elliptic curves $E_F$ over   number fields $F$ of small degree and discriminant
that are admissible for the Kummer construction, and for which our results on simultaneous resolutions apply. 
Ogg \cite{Ogg 1966} already  classified the elliptic curves over $F=\QQ$ where bad reduction occurs only at $p=2$. 
Over more general number fields,  such questions
attracted a lot of attention, for example in the work of Setzer \cite{Setzer 1978}, Stroecker \cite{Stroecker 1983}, 
Rohrlich \cite{Rohrlich 1982},  Ishii \cite{Ishii 1986}, 
Pinch \cite{Pinch 1984}, \cite{Pinch 1986}, \cite{Pinch 1998}, Bertolini and Canuto \cite{Bertolini; Canuto 1988},
Kida and Kagawa \cite{Kagawa 1998}, \cite{Kida; Kagawa 1997}, \cite{Kida 1999}, \cite{Kagawa 2001}, 
\cite{Kida 2001},  Cremona and Lingham \cite{Cremona; Lingham 2007}, Clemm and Trebat-Leder \cite{Clemm; Trebat-Leder 2016},
and myself \cite{Schroeer 2022}. In this regard, however,  little  is known about general curves, surfaces or higher-dimensional schemes.
 
The work of Pinch \cite{Pinch 1984} reveals that over $L=\QQ(\sqrt{-3})=\QQ(e^{2\pi i/3})$ there
is precisely one conjugate pair of elliptic curves $E_L$ that are admissible for the Kummer construction.
It turns out that this indeed leads to families $X\ra \Spec(\O_F)$ of K3 surfaces stemming from   $S_3$-number fields $F$.
The work of Comalada \cite{Comalada 1990} tells us that up to isomorphism,  there are exactly eight elliptic curves $E_F$
over quadratic number fields $F$ that have good reduction everywhere,  and the occurring
discriminants of the number fields are $d_F\in\{28,41,65\}$. 
Since the 2-torsion in $E_F$ is constant, no further base-change is needed to construct $X\ra\Spec(\O_F)$. 
The combination of  all these results yields the proof for Theorem A.

\begin{acknowledgement} 
I wish to thank the referees for many valuable suggestions, which helped to improve the paper.
The research was supported by the Deutsche Forschungsgemeinschaft in the project \emph{Varieties with Free Tangent Sheaf}, project
number 536205323. It was also conducted       in the framework of the   research training group
\emph{GRK 2240: Algebro-Geometric Methods in Algebra, Arithmetic and Topology}.
\end{acknowledgement}

\section{Simultaneous resolutions and \texorpdfstring{$\ell$}{l}-adic cohomology}
\mylabel{Simultaneous resolutions}

In this section we discuss the fundamental notations and results pertaining to simultaneous resolutions
of surface singularities, and the role of $\ell$-adic cohomology in this regard.

Let $S$ be a base scheme. So from now on, all schemes are endowed with a structure morphism to $S$,
and likewise all rings $R$ come with a   $\Spec(R)\ra S$, usually suppressed from notation.
A \emph{family of   normal surfaces} is a triple $(R,V,h)$ where $R$ is a ring,
$V$ is an algebraic space, $h:V\ra \Spec(R)$ is a proper flat morphism   of finite presentation
such that  for each geometric point
$s:\Spec(k)\ra \Spec(R)$, the base change $V_s=V\otimes_Rk$ is integral, normal, and two-dimensional.
Such families form a category fibered in groupoids over $(\Aff/S)$,
with the obvious notion of morphisms, and the stack axioms with respect to the fppf topology hold.
Also note that the    notion of families  
immediately generalizes to the situation where the ring $R$ is replaced by an arbitrary scheme
or algebraic space.

Let $V$ be a family of normal surfaces over $S$. We consider the category
$$
\uRes_{V/S}= \{(R,X,f,r)\}
$$
where $(R,X,f)$ is a family of  normal surfaces and 
$r:X\ra V\otimes R$ is an $R$-morphism  such that the following holds: For each geometric point $s:\Spec(k)\ra\Spec(R)$,
the fiber $X_s$ is regular, the morphism $r_s:X_s\ra V_s$ is birational, and the dualizing sheaf satisfies
$(\omega_{X_s}\cdot E)\geq 0$ for each irreducible curve $E$ in the exceptional locus $\Exc(X_s/V_s)$.
Then every other resolution  of singularities for $V_s$ uniquely factors over $X_s$.
Therefore, the objects are called
\emph{simultaneous minimal resolutions of singularities}.  

Note that for each point $s_0\in \Spec(R)$, the schematic fiber $V_0=f^{-1}(s_0)$ has  a minimal resolution of singularities $V'_0$,
which is  \emph{regular} in the sense that each local ring is regular, and unique up to isomorphism.
However, if $k_0=\kappa(s_0)$ is imperfect, $V'_0$ may fail to be \emph{geometrically regular}.
In this case,   $V'_0\ra\Spec(k_0)$ is not \emph{smooth}, and  
$V_0\ra\Spec(k_0)$ viewed as a family does not admit a simultaneous minimal resolution of singularities.
To cope with such situations, we follow
\cite{Fanelli; Schroeer 2020}, Section 1 and write $\Sing(V/S)$ for the closed subscheme defined by the second Fitting ideal
of $\Omega^1_{V/S}$, and call it the \emph{scheme of non-smoothness}. Its complement is denoted by $\Reg(V/S)$.
 
Clearly, the forgetful functor
$$
\uRes_{V/S}\lra(\Aff/S),\quad (R,X,f,r)\longmapsto \Spec(R)
$$  
is a category  fibered in groupoids, and satisfies the stack axioms. The fiber categories $\uRes_{V/S}(R)$ are actually
\emph{setoids}, that is, every arrow is an isomorphism and every self-arrow is an identity.
Indeed,  $\uRes_{V/S}$ is equivalent to the comma category of affine schemes
over an \emph{algebraic space} $\Res_{V/S}$, according to Artin's fundamental insight (\cite{Artin 1974}, Theorem 1),
which builds on Brieskorn's results (\cite{Brieskorn 1966}, \cite{Brieskorn 1968}, \cite{Brieskorn 1971}).
In particular, there is a \emph{universal simultaneous minimal resolution of singularities}
$$
r_\univ:X_\univ\lra \Res_{V/S}.
$$
Usually  the base, and thus also the total space, arehighly \emph{non-separated in an  un-schematic way}, compare the illustration 
in \cite{Artin 1974}, Figure 1.1.
From another perspective   the  algebraic space $\Res_{V/S}$ is rather close to being a scheme:

\begin{lemma}
\mylabel{purely inseparable extensions}
For each $s\in S$, the fiber $(\Res_{V/S})_s$ is the spectrum of a finite purely inseparable field
extension of the residue field $\kappa(s)$.
\end{lemma}

\proof
Without loss of generality we may assume $S=\Spec(k)$, and have to verify that $\Res_{V/k}$ is the spectrum of
some finite purely inseparable   extension $k'$. We start with some preliminary observations:

Let $k\subset L$ be some field extension, and write $r:X\ra V\otimes L$ for the minimal resolution of singularities.
As discussed above,  the scheme $X$ is regular, but not necessarily geometrically regular. Suppose the latter holds.
The closed   subscheme $Z=\Sing(X/V\otimes L)$ is one-dimensional,
but may have embedded components. Write $E$ for the closed subscheme corresponding to the sheaf of ideals $\shI\subset\O_Z$
comprising the torsion sections. Let $\eta_1,\ldots,\eta_s\in E$ be the generic points and $E_i$ be the schematic
closure of the monomorphism $\Spec(\O_{E,\eta_i})\ra X$. The intersection matrix $(E_i\cdot E_j)$ is negative-definite,
and we find some $D=-\sum \lambda_iE_i$  such that  $(D\cdot E_j)>0$. Fix a numerically ample Weil divisor $C\subset V$.
For some sufficiently large $\mu\geq 0$, we may assume that  Mumford's rational pullback $\mu r^*(C)$ has integral coefficients
(see \cite{Mumford 1961} and \cite{Schroeer 2019} for details), and furthermore $\mu r^*(C)+D$ is  ample on $X$. 
Let $\shL$ be the resulting ample sheaf. Passing
to some tensor power, it becomes very ample, and defines a closed embedding $X\subset\PP^m_L$. 
Let $k\subset k'\subset L$ be the common field of definition for the closed subscheme $X$ and the graph $\Gamma_r$
of the morphism $r:X\ra V\otimes L$.

We claim that $k\subset k'$ is finite and purely inseparable. For this we may enlarge $L$, and assume that $L$ is 
algebraically closed.   
Each member of  $\Gamma=\Aut(L/k)$ stabilizes $\Reg(V/k)\otimes L$, and therefore also its closure 
$X$ inside $\PP^m_L$. Thus $k'$ belongs to the fixed field $L^\Gamma$.
If $t\in L$ is transcendental, the automorphism group of $k(t)$ equals $\PGL_2(k)$, and each element extends
an  element of $\Gamma$ by \cite{Charnow 1970}, Lemma 1.
Likewise, for each separable $\omega\in L$
the resulting finite Galois extension has a non-trivial automorphism  group, and its members
can be extended to $\Gamma$. It follows that $L^\Gamma$ and hence $k'$ is purely inseparable.
Being a field of definition, $k'$ is finitely generated, and thus finite.

Suppose now that $k\subset L$ is purely inseparable. 
We claim that the extension $k\subset k'$ does not depend on $L$, up to unique isomorphism.
To see this, fix a perfect closure $k^\perf $, and  regard $L$ as a subfield of $k^\perf$.
The formation of $Z=\Sing(X/V\otimes L)$  commutes with all field extensions.
It follows from \cite{EGA IVb}, Proposition 4.2.7 that the same holds for the formation of the curves $E$.
Also the formation of $E_i$ does not change under purely inseparable extension.
Thus    $k'$ is the same, whether formed with  $L$ or $L'=k^\perf$.

We finally check that the  minimal resolution of singularities 
$r':X'\ra V\otimes k'$  is indeed universal, and thus $\Res_{V/k}=\Spec(k')$:
The preceding arguments already reveal that $|\Res_{V/k}|$ is a singleton.
The remaining task boils down to checking that  each simultaneous minimal resolution $X\ra V\otimes R$ with an 
Artin local ring $R$ is constant. This easily follows, by viewing  $X$ as a deformation of $X_0=X\otimes R/\maxid_R$
containing the constant family $\Reg(V/k)\otimes R$.
\qed

\medskip
We   may regard the topological space $|\Res_{V/S}|$, which comprises the equivalence classes
of morphisms $\Spec(K)\ra \Res_{V/S}$, as the set $|S|$, but  endowed with a possibly finer topology.
It is convenient to write $s'\in |\Res_{V/S}|$ for the point corresponding to $s\in S$,
and $\kappa(s')$ for the purely inseparable extension in Lemma \ref{purely inseparable extensions}, 
and   
\begin{equation}
\label{fiber universal resolution}
r_{s'}:X_{s'}=X_\univ\times_{\Res_{V/S}}\Spec\kappa(s') \lra V\times_S\Spec\kappa(s')=  V_{s'}
\end{equation} 
for the resulting minimal resolution of singularities.
The following two observations are   useful:
 
\begin{lemma}
\mylabel{rdp singularities}
Notation as above. If  $R^1r_{s',*}(\O_{X_{s'}})=0$ and   $V_{s'}$ is Gorenstein, then $R^1h_*(\O_{\tilde{W}})=0$
and $W=V_s$ is Gorenstein. Furthermore,   the minimal resolution \eqref{fiber universal resolution}   factors over  the base-change  
$\tilde{W}\otimes\kappa(s')$,
and the latter is normal.
\end{lemma}

\proof
The canonical map $\omega_{W\otimes\kappa(s')}\ra\omega_W\otimes\kappa(s')$ is bijective, so the dualizing sheaf  $\omega_W$ is invertible
and  the scheme  $W$ is Gorenstein.
Let $E_1,\ldots,E_r\subset\tilde{W}$ be the exceptional divisors, and write $K_{\tilde{W}/W}=\sum\lambda_i E_i$.
Since the resolution is minimal, we have $K_{\tilde{W}/W}\cdot E_i\geq 0$  and consequently $\lambda_i\leq 0$ for all indices
$1\leq i\leq r$, confer \cite{Schroeer 2019}, Proposition 2.4.
Seeking a contradiction, we assume that  $R^1h_*(\O_{\tilde{W}})$ is non-zero. By \cite{Artin 1966}, Theorem 3
the fundamental cycle $Z\subset \tilde{W} $ has $\chi(\O_Z)\leq 0$. From 
$$
0\leq -2\chi(\O_Z)=\deg(\omega_Z) = (K_{\tilde{W}}+Z)\cdot Z<K_{\tilde{W}/W}\cdot Z
$$ 
we see that $\lambda_j<0$ for at least one $1\leq j\leq r$. In turn,
the relative canonical divisor for the normalization of  $\tilde{W}\otimes\kappa(s')$ with respect to $W\otimes\kappa(s)$
likewise has a strictly negative coefficient. This normalization is dominated by some regular modification of $X_{s'}$.
By our assumptions we have $K_{X_{s'}/V_{s'}}=0$, hence on the regular modification all coefficients are $\geq 0$,
contradiction.

This establishes $R^1h_*(\O_{\tilde{W}})=0$. We infer that all coefficients $\lambda_i$ vanish, 
and also    that the base-change  $\tilde{W}\otimes\kappa(s')$ is normal. If the latter is dominated by
some modification   rather than the minimal  resolution $X_{s'}$, then $\lambda_j>0$ for at least one $1\leq j\leq r$,
contradiction.
\qed

\medskip
Sometimes the condition on the singularities is enforced by a global hypothesis:

\begin{lemma}
\mylabel{rdp singularities by global condition}
Notation as above. If the dualizing sheaf $\omega_{\tilde{W}}$ is anti-nef, we have   
$R^1h_*(\O_{\tilde{W}})=0$. If $\omega_{\tilde{W}}$ is  numerically trivial, $W=V_s$ is furthermore Gorenstein.
\end{lemma}

\proof
Let $Z\subset W$ be the fundamental cycle. Then 
$-2\chi(\O_Z) =  (K_{\tilde{W}}\cdot Z)+ Z^2 < 0$,
so \cite{Artin 1966}, Theorem 3 ensures that $R^1h_*(\O_{\tilde{W}})=0$. Suppose now that
$K_{\tilde{W}}= r^*K_W + K_{\tilde{W}/W}$ is numerically trivial, and let $E_1,\ldots, E_r\subset W$ be the
irreducible components. Then $(K_{\tilde{W}/W}\cdot E_i)=0$, hence  $\omega_{\tilde{W}}$ is numerically trivial
on $Z$. 
Using $R^1h_*(\O_{\tilde{W}})=0$ we see that it is actually trivial on $Z$ and all its infinitesimal neighborhoods,
and infer that $\omega_W$ is invertible.
\qed

\medskip
Recall that on normal surfaces, the \emph{rational Gorenstein singularities}  are precisely the \emph{rational double points}.
Let us say that $V\ra S$ is an \emph{RDP family of normal surfaces} if all geometric fibers
$V_s$ are normal surfaces with at most rational double points.
We next examine the local henselian ring $\O^h_{\Res_{V/S},s'}$ whose residue field
$\kappa(s')$ is a finite purely inseparable extension of $\kappa(s)$. The ring comes  with a homomorphism
$$
\O_{S,s}\subset \O_{S,s}^h\lra \O^h_{\Res_{V/S},s'}.
$$
In general, this arrow is neither injective nor surjective, and it could well be that
the local ring for $s'$ coincides with the residue field. Now Artin's reformulation 
(\cite{Artin 1974}, Theorem 2) of Brieskorn's result is:

\begin{theorem}
\mylabel{artin-brieskorn}
If   $V\ra S$ is an RDP family of normal surfaces, the induced maps
$\Spec(\O^h_{\Res_{V/S},s'})\ra \Spec(\O_{S,s}^h)$ are surjective.  
\end{theorem}
 
In what follows, we freely use the theory of $\ell$-adic cohomology, also in the context of algebraic spaces.
\emph{To simplify the exposition, we now  assume  that $V\ra S$ is an RDP family of normal surfaces,
and that $S$ is the spectrum of a henselian discrete valuation ring $R$, with generic point $\eta$ and closed point $\sigma$.}

Note that $|S|=\{\eta,\sigma\}$ and  $|\Res_{V/S}|=\{\eta',\sigma'\}$,   that the canonical map between
these topological spaces is a homeomorphism, and  that  there is a universal simultaneous minimal resolution $f_\univ:X_\univ\ra\Res_{V/S}$.
We now fix a prime $\ell>0$ that is invertible in $R$, and consider the ensuing $\ell$-adic sheaf
$R^2f_{\univ,*}(\ZZ_\ell(1))$ on the algebraic space $\Res_{V/S}$.
The cycle classes for the exceptional divisors on $X_{\sigma'}\otimes \kappa(\sigma')^\sep$ 
define  an $\ell$-adic subsheaf inside $\uH^2(X_{\sigma'},\ZZ_\ell(1))$. The latter denotes the
$\ell$-adic sheaf on the spectrum of $\kappa(s')$.  By Lemma \ref{extensions local systems} below,
the inclusion  extends to an  $\ell$-adic subsheaf
$$
\Psi_{V/S}\subset R^2f_{\univ,*}(\ZZ_\ell(1)).
$$
This is an $\ell$-adic local system whose rank $r\geq 0$ is given by the number of exceptional divisors
in $X_{\sigma'}\otimes \kappa(\sigma')^\sep$.
If the generic fiber of $V$ is already smooth,  restriction to $\eta\in S$ gives
an $\ell$-adic subsheaf 
$$
(\Psi_{V/S})_{\eta} \subset\uH^2(X_{\eta'},\ZZ_\ell(1)) = \uH^2(V_{\eta},\ZZ_\ell(1))
$$
over the field of fractions $F=\Frac(R)$.
A powerful result of Shepherd-Barron (\cite{Shepherd-Barron 2021}, Corollary 2.14) now can be reformulated in the following way:

\begin{theorem}
\mylabel{shepherd-barron}
In the above situation, suppose the following conditions:
\begin{enumerate}
\item 
The generic fiber $V_\eta$ is smooth.
\item
The minimal resolution of singularities for $V_\sigma$ is smooth. 
\item
All exceptional divisors for $r_\sigma:X_\sigma\ra V_\sigma$ are isomorphic to $\PP^1_\sigma$.
\item
The local system $(\Psi_{V/S})_\eta$ is constant.
\end{enumerate}
Then the family $V\ra S$ admits a simultaneous resolution of singularities.
\end{theorem}

Note that for the assumptions and the conclusion alike there is no need to mention the algebraic space $\Res_{V/S}$.
To apply the result, the crucial point is to identify  the local system $(\Psi_{V/S})_\eta$ in condition (iv) from
the geometry of the family $V\ra S$. Also note that condition (ii) 
simply means that the purely inseparable extension $\kappa(\sigma)\subset \kappa(\sigma')$ is an equality.

Let $Y$ be a scheme or more generally an algebraic space that is connected and non-empty, $y_0:\Spec(\Omega)\ra Y$ be a morphism with some separably closed field $\Omega$,
and $\ell>0$ be a prime that is invertible on $Y$. Let 
$F=(F_\nu)$ be an \emph{$\ell$-adic local system}. The entries are \emph{$\ell^n$-local systems} of some rank $r\geq 0$,
in other words   abelian sheaves on the \'etale site $(\Et/Y)$
that are twisted forms of   the constant sheaves $(\ZZ/\ell^\nu\ZZ)_Y^{\oplus r}$.
Write $M=\invlim F_\nu(y_0)$ for the corresponding   $\ell$-adic representations   
of the fundamental group $\pi_1(Y,y_0)$.  Suppose there is a field $K$ such that $y_0$ factors over some monomorphism $\Spec(K)\ra Y$,
and write $ K^\sep \subset\Omega$ for the relative separable closure. Composing with the canonical map 
\begin{equation}
\label{galois to fundamental}
\Gal(K^\sep/K)^\op=\pi_1(\Spec K,y_0)\lra \pi_1(Y,y_0),
\end{equation} 
our representation  of the fundamental group $M$ becomes a Galois $\ZZ_\ell$-module.

\begin{lemma}
\mylabel{extensions local systems}
In the above situation, let  $M'\subset M$ be   Galois $\ZZ_\ell$-submodule.
If the canonical map \eqref{galois to fundamental} is surjective, then there is a unique $\ell$-adic subsystem  $F'\subset F$
with $M'=\invlim F'_\nu(y_0)$.
\end{lemma}

\proof
The assumption ensures that the $\ZZ_\ell$-submodule $M'\subset M$ is stable by the action of the fundamental group. 
Being a subrepresentation with respect to $\pi_1(Y,y_0)$, it   stems from a local subsystem,
because $F\mapsto \invlim F_\nu(y_0)$ is an equivalence of categories.
\qed

\medskip
%
We will apply the above results with  \emph{K3 surfaces}. 
Recall that over an algebraically closed ground field $k$ of characteristic $p\geq 0$, these are 
the smooth proper surfaces $X$ with $h^0(\O_X)=1$ and $c_1=0$ and $b_2=22$ 
(\cite{Bombieri; Mumford 1977} and \cite{Bombieri; Mumford 1976}).
An important variant: A \emph{normal K3 surface} is a normal proper surface $Y$ such that  the minimal resolution of singularities $r:X\ra Y$ 
is a K3 surface.  Using Lemma \ref{rdp singularities by global condition} we then get
$$
h^0(\O_Y)=h^2(\O_Y)=1,\,  h^1(\O_Y)=0\quadand  R^1r_*(\O_X)=0 \quadand \omega_Y=\O_Y.
$$
In particular, the exceptional divisors $E_i\subset X$, $1\leq i\leq r$ have dual graph whose connected components correspond
to the Dynkin diagram $A_n$ or $D_n$ or $E_6$, $E_7$, $E_8$.

We now replace the ground field $k$ by a base scheme $S$, and examine the  ensuing relative notions: A \emph{family of K3 surfaces} over  $S$ 
is an algebraic space $X$, together with a proper flat morphism  $f:X\ra S$   of finite presentation,
such that for all geometric points $s:\Spec(k)\ra S$ the fiber $X_s$ is a K3 surface  in the sense of the preceding paragraph.
The \emph{stack of K3 surfaces} $\shM_\Kthree\ra(\Aff/S)$ is formed by the triples $(R,X,f)$ where $R$ is a ring endowed with a structure
morphism $\Spec(R)\ra S$, and $f:X\ra \Spec(R)$ is a family of K3 surfaces.
Likewise,  one can consider  \emph{families of normal K3 surfaces}.  Note that these are automatically RDP families.
Furthermore, the notions immediately extend when we replace the base scheme $S$ by an algebraic space.

\section{Effective models for sign involutions}
\mylabel{Effective models}

Let $R$ be a discrete valuation ring, with maximal ideal $\maxid$, residue field $k=R/\maxid$,   
and field of fractions $F=\Frac(R)$. 
Write $p\geq 0$ for the characteristic of the residue field, and choose a uniformizer $\pi\in R$.
Let 
\begin{equation}
\label{weierstrass equation}
y^2+a_1xy+a_3 y = x^3+a_2x^2 + a_4x + a_6
\end{equation} 
be a Weierstra\ss{} equation with coefficients $a_i\in R$.  
One also has   auxiliary  values and discriminant $b_i,c_i,\Delta\in R$, as defined in \cite{Deligne 1975}, Section 1.
The homogeneous cubic equation
\begin{equation}
\label{homogeneous weierstrass equation}
Y^2Z+a_1XYZ + a_3YZ^2 = X^3+a_2X^2Z +   a_4XZ^2+a_6Z^3
\end{equation} 
defines a family of cubic curves $E\subset \PP^2_R$ containing
the  \emph{zero section}
$$
\Zero(E/R) = V_+(X,Z) =\{(0:1:0)_R\}.
$$ 
The dehomogenization via $x=X/Z$ and $y=Y/Z$
gives back our Weierstra\ss{} equation, which defines the affine part $E\cap D_+(Z)$.
As customary, we write $D_+(\ideala)$    for the open set defined by a homogeneous ideal, 
and $V_+(\ideala)$ for the complementary closed set. 

Following \cite{Fanelli; Schroeer 2020}, Section 1, we write $\Sing(E/R)$ for the closed subscheme defined by the first Fitting ideal
of $\Omega^1_{E/R}$, and call it the \emph{scheme of non-smoothness}. The  complementary open subscheme
$\Reg(E/R)$ is the smooth locus, and carries a   group law in  a canonical way. Its fibers are either elliptic curves, or twisted forms of the 
multiplicative group $\GG_m$,  or twisted forms of the additive group $\GG_a$,
in dependence on  $\Delta,c_4\in R$.  Note that we have 
\begin{equation}
\label{equivalent vanishing values}
2,a_1,a_3 =0\quad\Longleftrightarrow\quad 2, \Delta,c_4=0,
\end{equation} 
which follows from the congruences $c_4\equiv a_1^4$ modulo 2,  together with $\Delta\equiv a_3^4$ modulo the ideal $(2,a_1)$.

\emph{Throughout, we assume that at least one of  the elements $2,a_1,a_3\in R$
is non-zero}, or equivalently one of $2,c_4,\Delta\in R$ is non-zero. Another way to phrase this is that the matrix
$$
N=
\begin{pmatrix}
1	& -a_1	& 0\\
0	& -1	& 0\\
0	& -a_3	& 1
\end{pmatrix}
\in\GL_3(R)
$$
has order two, and   the same holds for the  image in $\PGL_3(R)=\Aut(\PP^2_R)$. The resulting involution of $\PP^2_R$ 
stabilizes the closed subscheme $E$,
and induces the negation map   on the commutative group scheme $\Reg(E/R)$. 
We therefore call   $N:E\ra E$ the   \emph{sign involution}.
On the affine part, it is given by $x\ra x$ and $y\ra -(y+a_1x+a_3)$.

According to \cite{FGA IV}, Section 4 the functor of relative automorphisms for the family $E\ra\Spec(R)$  
is representable by
a relative group scheme $\Aut_{E/R}$ that is separated and  of finite type, and the sign involution defines a homomorphism
$(\ZZ/2\ZZ\cdot N)_R\ra \Aut_{E/R}$. The following terminology follows Romagny (\cite{Romagny 2012}, Definition 4.3.3):

\begin{definition}
\mylabel{effective model}
The \emph{effective model} for the sign involution is the schematic closure
$$
G=\overline{(\ZZ/2\ZZ\cdot N)_F}\subset\Aut_{E/R}
$$
of the generic sign involution inside the relative automorphism group scheme.
\end{definition}

According to loc.\ cit., Theorem 4.3.4 the structure morphism $G\ra\Spec(R)$ is finite and flat, of degree two,
and the inclusion  $G\subset\Aut_{E/R}$ is a closed subgroup scheme.

Recall that over any ring $A$, given two integers $a,b\in A$ with $ab=2$ we get the group valued functor
$$
G^{a,b}(A')=\{f\in A'\mid f^2-af=0\}
$$
with group law  $f_1\ast f_2=f_1+f_2- bf_1f_2$. 
According to the Tate--Oort Classification (\cite{Tate; Oort 1970}, Theorem 2), any family of group schemes over $A$
of order two is isomorphic to some $G^{a,b}$.
Moreover, we have $G^{a,b}\simeq G^{a',b'}$ if and only if
$(a,b)$ and $(a',b')$ lie in the same orbit of the action of $A^\times$ on 
$A\times A$  given by $\lambda\cdot (u,v)=(\lambda u, \lambda^{-1}v)$.
Note that the pairs  $(a,b)=(1,2)$ and $(a,b)=(2,1)$ define the constant group
scheme $(\ZZ/2\ZZ)_A$ and the multiplicative group scheme $\mu_{2,A}$, respectively.
If $2=0$ in $A$, we have a further pair $(a,b)=(0,0)$, which yields the unipotent group scheme $\alpha_2$.
Throughout, we write 
$$
C=\Spec k[u^2,u^3]\cup\Spec k[u^{-1}]
$$ 
for the \emph{rational cuspidal curve} over the field $k$.

\begin{proposition}
\mylabel{infinitesimal specialization}
The canonical map $(\ZZ/2\ZZ\cdot N)_R\ra G$ is not an isomorphism if and only if   $2,a_1,a_3\in\maxid$.
In this case, the residue field $k$ has characteristic two, 
and the closed fiber $G_k$ of the effective model is a copy of  $\mu_{2,k}$ or $\alpha_{2,k}$,
and the genus-one curve $E_k$ is a twisted form of $C$.
\end{proposition}

\proof
Clearly, the map is not an isomorphism if and only if the sign involution on the affine part of $E_k$ is trivial,
in other words, if $-(y+a_1x+a_3)=y$. This means $-1\equiv 1$ and $-a_1\equiv 0$ and $-a_3\equiv 0$ modulo $\maxid$,
and the first assertion follows. 

Suppose these equivalent conditions hold. Then   the residue field $k$ has characteristic two.
Since the generic fiber  $(\ZZ/2\ZZ)_F$ is dense in both $(\ZZ/2\ZZ)_R$ and $G$,
it follows that $G_k$ is infinitesimal. Up to isomorphism, the only infinitesimal group schemes of order two are  $\mu_{2,k}$ or $\alpha_{2,k}$.
From \eqref{equivalent vanishing values} we get $ \Delta,c_4\in \maxid_R$. So by \cite{Deligne 1975}, Proposition 5.1 the closed fiber $E_k$
is a twisted form of the rational cuspidal curve.
\qed

\medskip
Recall that the \emph{scheme of $N$-fixed points} is the intersection of graph  $\Gamma_N$
with the diagonal $\Delta_{\PP^2_R}$ inside $\PP^2_R\times \PP^2_R$.
Its formation commutes with base-change in $R$.

\begin{proposition}
\mylabel{fixed scheme}
The scheme of $N$-fixed points $\Gamma_N\cap\Delta_{\PP^2_R}$ is defined inside the diagonal by the homogeneous ideal
$\ideala=(X,Z)\cdot (2Y+a_1X+a_3Z)$. Moreover, the generic fiber of $V_+(\ideala)\cap E$ is finite of degree four.
\end{proposition}

\proof
First consider the invariant open set $D_+(Z)$. With respect to $x=X/Z$ and $y=Y/Z$,
the action is given by $x\mapsto x$ and $y\mapsto -y-a_1x-a_3$. The ideal for the scheme
of fixed points is defined by $2y+a_1x+a_3=0$. In the same way one sees that on $D_+(X)$
the scheme of fixed points is defined by  $2y'+a_1+a_3z'=0$.
It follows that  $\Gamma_N\cap\Delta=V_+(\ideala)$  holds  over  the open set $D_+(Z)\cup D_+(X)$. 

The complementary closed set is  $\Zero(E/R)$. Let $\sigma$ be its closed point, 
and write $A=\widehat{R}[[X/Y,Z/Y]]$ for the formal completion of the local ring $\O_{\PP^2_R,\sigma}$.
Then the  monomorphism $\Spec(A)\ra\PP^2_R$ is stable with respect to $N$,
and the matrix acts via 
$$
\frac{X}{Y}\longmapsto  \frac{X}{Y}\cdot\frac{-1}{1+a_1X/Y + a_3Z/Y}\quadand \frac{Z}{Y}\longmapsto  \frac{Z}{Y}\cdot\frac{-1}{1+a_1X/Y + a_3Z/Y}.
$$
Now the  scheme of fixed points is defined by $-1=1+a_1X/Y + a_3Z/Y$.
This shows $\Gamma_N\cap\Delta=V_+(\ideala)$ after base-changing with respect to $\Spec(A)\ra\PP^2_R$.
The disjoint union $D_+(Z)\cup D_+(X)\cup \Spec(A)$ is an fpqc covering, and our assertion follows.

It remains to check the assertion on the generic fiber. If $2\neq 0$ in $R$, then $V_+(\ideala)\otimes F$
is the disjoint union of the point $(0:1:0)$ and the line $V_+(2Y+a_1X+a_3Z)$. Using Bezout's Theorem,
we infer that the intersection with $E_F$ is finite of degree $4=1+3$.
If $\Delta\neq 0$ the intersection   is the 2-torsion part of the elliptic curve $E_F$, which has degree $4=2^2$.
Finally, suppose that both $2,\Delta\in R$ vanish. 
Using \eqref{equivalent vanishing values} we get $c_4\neq 0$. From  \cite{Deligne 1975}, Proposition 5.1 we see that $\Reg(E_F/F)$ is a twisted form of $\GG_{m,F}$.
The two-torsion part $\mu_2$ has degree two. The formal completion at the singular point in $E_F$
is a twisted form of $F[[u,v]]/(uv)$, and   sign involution permutes the indeterminates.
The fixed scheme is given by the relations $uv=0$ and $u=v$, which also has degree two.
Again our fixed scheme has degree $4=2+2$.
\qed

\medskip
Recall that   \emph{schemes of fixed points} exist  in surprising generality (\cite{Demazure; Gabriel 1970}, Chapter II, \S1, Theorem 3.6).
In particular, we can form the  scheme of  fixed points
$$
\Fix(E/R)=E^G
$$ 
for the effective model $G$ of the sign involution on the family of cubics $E$.
It commutes with base-change, and  is related to the scheme of $N$-fixed points, which   may contain the closed fiber and therefore fails to be flat, as follows:
Inside the principal ideal domain $R$ one has  $\gcd(2,a_1,a_3) = \pi^d$ for some unique $d\geq 0$. As in \cite{Tate 1975}, page 48   we write 
$$
2_d=2/\pi^d\quadand a_{1,d}=a_1/\pi^d\quadand a_{3,d}=a_3/\pi^d,
$$
and introduce  the   homogeneous ideal $\idealb=(X,Z)\cdot (2_dY + a_{1,d}X + a_{3,d}Z)$.

\begin{proposition}
\mylabel{effective fixed scheme}
We have $\Fix(E/R)=E\cap V_+(\idealb)$. Moreover, this scheme is finite and   flat over $R$, of relative degree four.
\end{proposition}

\proof
According to Proposition \ref{fixed scheme}, the statements hold  after base-changing along $R\subset F$.
We next examine  the closed fiber of $B=E\cap V_+(\idealb)$.
If $2_d\not\in\maxid$ we can argue as in the proof for Proposition \ref{fixed scheme}.
Suppose now $2_d$ vanishes in $k=R/\maxid$. Then $V_+(\idealb)_k$
contains the line $L=V_+(a_{1,d}X + a_{3,d}Z)_k$, and the corresponding sheaf of ideals
is the torsion sheaf   supported by $(0:1:0)$ with stalk $k$, so its degree is $l=1$. Using 
Bezout's Theorem, we see that $B_k$ is finite, of degree at most  $3+l=4$.
Since $B\ra\Spec(R)$ is proper with finite fibers, it must be finite. Its coordinate
ring $A=\Gamma(B,\O_B)$ is a finitely generated $R$-module of rank four
with $\dim_k(A\otimes_Rk)\leq 4$. Using the structure theory for such modules,
we see that $A$ is free.

Consider the morphism $h:G\times B\ra E$ stemming from the $G$-action on $E$.
Since $(G\times B)_F$ is schematically dense in $G\times B$,
the schematic image $h(G\times B)\subset E $ is the closure of $B_F$. 
This equals $B$, by flatness of $B$, and we infer that $B\subset E$ is $G$-stable.
By construction, the $G$-action is trivial on $B_F$.  We thus have $h_F=(\pr_2)_F$.
Using schematic density of $(G\times B)_F$ again, we see $h=\pr_2$, thus the $G$-action on $B$
is trivial. 

Summing up, we have an inclusion $B\subset \Fix(E/R)$, which is generically an equality. It remains
to check that the closed fiber of the scheme of fixed points has degree four.
This follows as in the proof for Proposition \ref{fixed scheme} if one of $\Delta,c_4$ does not belong to $\maxid$.
Otherwise, $E_k$ is a twisted form of the rational cuspidal curve
$C$, and it suffices to treat the case $E_k=C$. 
If $2\not\in\maxid$ the sign involution is given by $u^{-1}\mapsto -u^{-1}$, and one immediately sees
that the fixed scheme is given on the two charts by the respective relations $u^3=0$ and $u^{-1}=0$,
so its length is $3+1=4$. The case $2\in\maxid$ follows from   Proposition \ref{fixed scheme determines action} below.
\qed

\medskip
Consider the rational cuspidal curve $C$ over a field $k$ of characteristic two.
Despite the singularity, the   tangent sheaf $\Theta_{C/k}$  is invertible,  having degree four (\cite{Schroeer 2007}, Section 3).
The scheme of fixed points already determines the closed fiber $G_k$ of the effective model, at least in the case 
that matters to us:

\begin{proposition}
\mylabel{fixed scheme determines action}
In the above setting, the  map $H\mapsto C^H$ gives a bijection between the infinitesimal subgroup schemes $H\subset\Aut_{C/k}$ of order two,
and the effective Cartier divisors $D\subset C$ of degree four with $\O_C(D)\simeq\Theta_{C/k}$. Moreover, we have $H\simeq\mu_2$ if and only
if $D$ is \'etale at some point.
\end{proposition}

\proof
The one-dimensional vector subspaces $L=k\delta$ inside $H^0(C,\Theta_{C/k})$ correspond to the effective Cartier divisors 
$D\subset C$   with $\O_C(D)\simeq\Theta_{C/k}$.
Every such $L$  is stable with respect to Lie bracket and $p$-map
(\cite{Schroeer 2007}, Section 3   and \cite{Kondo; Schroeer 2021}, Section 2).  
So by the Demazure--Gabriel Correspondence (\cite{Demazure; Gabriel 1970}, Chapter II, \S7, No.\ 4),  
the $L$ also correspond to the infinitesimal subgroup schemes $H$ of order two.
Moreover, the scheme of fixed points coincides with the zero scheme for $\delta\in H^0(C,\Theta_{C/k})$, in other words $D=C^H$.

As explained in \cite{Schroeer 2007}, proof of Proposition 3.2 the derivations
$u^{-2}D_u,  uD_u,  D_u,  u^2D_u$ form a basis for $\lieg=H^0(C,\Theta_{C/k})$. The discussion in \cite{Kondo; Schroeer 2021}, Section 2
reveals that $H\simeq\mu_2$ if and only if $C^H$ is \'etale.
\qed

\section{Freeness at the singular locus}
\mylabel{Freeness}

Keep the setting of the previous section, such that $E$ is a family of cubic curves inside  $ \PP^2_R=\Proj R[X,Y,Z]$,
with at least one of the $2,a_1,a_3\in R$ non-zero, 
and $G$ is the effective model of the sign involution.
To simplify exposition, we also assume that the generic fiber $E_F$ is smooth, in other words, the discriminant $\Delta\in R$ is non-zero.
The zero scheme, the scheme of non-smoothness, and the scheme of fixed points 
$$
\Zero(E/R)\quadand \Sing(E/R)\quadand \Fix(E/R),
$$
define closed subschemes of $E$, 
all of which are  finite over $R$.   The goal of this section is to clarify how these three schemes intersect.
Throughout, $\val:R\ra\ZZ\cup\{\infty\}$ denotes the valuation, normalized by $\val(\pi)=1$.
Recall that we write $2_d=2/\pi^d$ and $a_{1,d}=a_1/\pi^d$ and $a_{3,d}=a_3/\pi^d$,   where $d=\val(2,a_1,a_3)$.
One says that the closed fiber $E_k$ has \emph{additive type} if it is a twisted form of  the rational cuspidal curve.

\begin{proposition}
\mylabel{fix and singular}
The scheme $\Fix(E/R)\cap\Sing(E/R)$ is disjoint from $\Zero(E/R)$, and is given as a closed subscheme in  
$\AA^2_R=D_+(Z)$ by the   equations
\begin{equation}
\label{fix and singular locus}
-y^2=x^3+a_2x^2 + a_4x + a_6,\quad
a_1y=3x^2+2a_2x+a_4,\quad 
2_dy + a_{1,d}x+   a_{3,d} =0.
\end{equation} 
If the residue field $k=R/\maxid$ has characteristic two and  the closed fiber $E_k$ has additive type,
the  intersection is set-theoretically described by
$$
\pi=0,\quad x^2=a_4,\quad y^2=a_2a_4+a_6,\quad 2_dy + a_{1,d}x+   a_{3,d} =0.
$$
\end{proposition}

\proof
Taking   partial derivatives of the homogeneous Weierstra\ss{} equation \eqref{homogeneous weierstrass equation}, one already  sees that  
$\Sing(E/R)$ is disjoint from  $\Zero(E/R)$. In turn, the intersection $\Fix(E/R)\cap\Sing(E/R)$ is contained
in $\AA^2_R=D_+(Z)$. Thus $\Sing(E/R)$ is defined by the Weierstra\ss{} equation $y^2+a_1xy+a_3 y = x^3+a_2x^2 + a_4x + a_6$,
together with the equations $a_1y = 3x^2+2a_2x+a_4$ and $ 2y+a_1x+a_3=0$
coming from the partial derivatives. The latter is a multiple of 
$2_dy + a_{1,d}x+   a_{3,d} =0$, which by Proposition \ref{effective fixed scheme} defines  $\Fix(E/R)\smallsetminus\Zero(E/R)$.
Thus we already have the second and third equations in \eqref{fix and singular locus}. 
The first is obtained from the Weierstra\ss{} equation by subtracting a multiple of the third.
 
For the remaining statement, suppose  that  $k$ has characteristic two and $E_k$ has additive type.
Since    $\Sing(E/R)$ is set-theoretically contained in the closed fiber,
we can enlarge \eqref{fix and singular locus} by the equation  $\pi=0$ and still describe  $\Fix(E/R)\cap\Sing(E/R)$
as a closed set.
Moreover, we have $\Delta,c_4\in\maxid$ according to  \cite{Deligne 1975}, Proposition 5.1. From \eqref{equivalent vanishing values}
we get $a_1,a_3\in\maxid$.
Using this for substitution in  \eqref{fix and singular locus},  the second equation  becomes $x^2=a_4$, and the first   $y^2=a_2a_4+a_6$.
\qed

\medskip
We now can determine when the intersection in Proposition \ref{fix and singular} is empty.
As in \cite{Tate 1975}, Section 7, it is convenient to write $\bar{a}_i, \bar{a}_{j,d}\in k$ et cetera for the residue classes
for $a_i,a_{j,d}\in R$.

\begin{proposition}
\mylabel{no intersection}
Suppose the closed fiber $E_k$ is non-smooth.
Then the  subschemes  $\Fix(E/R)$ and $\Sing(E/R)$ are disjoint if and only if the following three conditions hold:
\begin{enumerate} 
\item The closed fiber $E_k$ is of additive type.
\item The residue field $k=R/\maxid$ has characteristic two.
\item We have $\val(b_2a_4+b_6)=2\val(2,a_1,a_3)$. 
\end{enumerate}
Furthermore, condition (iii) is equivalent to $\val(a_2a_4+a_6)=0$, provided  $a_1=a_3=0$.
\end{proposition}

\proof
We start with some preliminary observations:
First note that $\Sing(E/R)=\{z\}$ is a singleton whose residue field extension $k\subset \kappa(z)$
is purely inseparable. It follows that  $z$ belongs to the fixed scheme of any automorphism.
Next note that from  the very definition of   $b$-values in \cite{Tate 1975}, Section 1  we    get the relation
\begin{equation}
\label{relation b-values a-values}
b_2a_4+b_6 =  2^2(a_2a_4+a_6) + a_1^2a_4 + a_3^2.
\end{equation} 
So if $a_1=a_3=0$ we have $\val(b_2a_4+b_6)=2\val(2,a_1,a_3)+\val(a_2a_4+a_6)$  and   see that (iii) simplifies as claimed.
Finally,   suppose that $k=R/\maxid_R$ has characteristic two, and let 
$B\subset\AA^2_k$ be the closed subscheme defined by $x^2=\bar{a}_4$ and 
$y^2=\bar{a}_2\bar{a}_4+\bar{a}_6$ and $2_dy + \bar{a}_{1,d}x+   \bar{a}_{3,d} =0$.
Substituting the former two in the square of the latter, we immediately see
$$
B=\varnothing\quad \Longleftrightarrow\quad(2_d)^2(\bar{a}_2\bar{a}_4+\bar{a}_6) +  \bar{a}_{1,d}^2\bar{a}_4 +  \bar{a}_{3,d}^2\neq 0.
$$
We now can easily check that our  conditions (i)--(iii)  are necessary: Suppose that $\Fix(E/R)$ and $\Sing(E/R)$ are disjoint.
If     the closed fiber is  of multiplicative type, or if  the residue field $k=R/\maxid$ has characteristic $p\neq 2$,
then  the effective model of the sign involution is $G=(\ZZ/2\ZZ\cdot N)_R$, and thus $z\in \Fix(E/R)$, contradiction.
This gives (i) and (ii). From Proposition  \ref{fix and singular}, it then follows   $B=\Fix(E/R)\cap\Sing(E/R)$.
This is empty, hence the fraction
$$
\frac{b_2a_4+b_6}{\pi^{2d}}=\frac{2^2(a_2a_4+a_6) + a_1^2a_4 + a_3^2}{\pi^{2d}} 
$$
belongs to $R^\times$, and consequently (iii) holds.
 
The conditions  (i)--(iii) are also  sufficient: The third gives $B=\varnothing$, by the observations
in the preceding paragraph. Conditions (i) and (ii) ensure that $B$ coincides with the intersection $\Fix(E/R)\cap \Sing(E/R)$,
according to  Proposition \ref{fix and singular}.
\qed

\medskip
Let us finally determine when the closed fiber of the effective model is a multiplicative group scheme:

\begin{proposition}
\mylabel{fixed scheme etale}
Suppose the closed fiber $E_k$ is non-smooth, and   $\Fix(E/R)$ and $\Sing(E/R)$ are disjoint.
Then the following are equivalent:
\begin{enumerate}
\item
The fixed scheme $\Fix(E/R)$ is \'etale over $R$.
\item
The group scheme $G_k$ is isomorphic to $\mu_{2,k}$.
\item
We have $\val(2)\leq \val(a_i)$ for $i=1$ and $i=3$.
\end{enumerate}
\end{proposition}

\proof
First observe that we  already saw (i)$\Leftrightarrow$(ii) in Proposition \ref{fixed scheme determines action},
and that condition (iii) is equivalent to $2_d\not \in \maxid$.
By Proposition \ref{effective fixed scheme}, the fixed scheme $\Fix(E/R)_k$ 
has multiplicity $m\geq 2$ at the origin $(0:1:0)$ if and only if $2_d\in\maxid$.
So   (iii)$\Leftrightarrow$(ii) follows from Proposition \ref{fixed scheme determines action}.
\qed

\section{Kummer constructions in families}
\mylabel{Kummer constructions}

Throughout this section,  $S$ denotes  an irreducible   Dedekind scheme  with generic point $\eta\in S$,
and $E_F$ is an  elliptic curve over the function field $F=\kappa(\eta)$.
Note that we allow imperfect residue fields $\kappa(s)$. 
To simplify exposition we assume that $S$ is \emph{excellent}, which here means that for 
each closed point  $s\in S$ the formal completion $\O_{S,s}\subset \O_{S,s}^\wedge$    induces a separable   extension on field of fractions 
(compare \cite{Schroeer 2020}, Proposition 4.1). This obviously holds if $F=\Frac(\O_{S,s})$ has characteristic zero, the case
we are mainly interested in.
The goal of this  section is to introduce the \emph{Kummer construction}
in a relative setting over $S$, which leads to families of normal K3 surfaces.

Let  $E\ra S$ be the  \emph{Weierstra\ss{} model} for the elliptic curve $E_F$. This is obtained
from the \emph{N\'eron model}  by passing to the minimal regular compactification, and contracting
all curves disjoint from the zero section. 
For each closed point $s\in S$ the base change $E\otimes\O_{S,s}$
becomes a family of cubic curves defined by a Weierstra\ss{} equation 
$$
y^2 + a_1xy + a_3y = x^3 + a_2x^2+a_4x+a_6
$$
over the discrete valuation ring $R=\O_{S,s}$.
We refer to  \cite{Tate 1975} or \cite{Silverman 1994} for details on the 
auxiliary values and discriminant $b_i, c_i, \Delta\in R$, the  \emph{invariant}   $j=c_4^3/\Delta$,
the valuations $v_s=\val_s(\Delta)\geq 0$ and the \emph{discriminant divisor} $\sum v_s\cdot s$,
which supports the \emph{points of bad reduction}.

By the universal property of N\'eron models (\cite{Bosch; Luetkebohmert; Raynaud 1990}, Section 1.2), the sign involution on the generic fiber $E_F$ 
extends to a homomorphism $(\ZZ/2\ZZ)_S\ra\Aut_{E/S}$. As in Section \ref{Effective models}, the schematic closure
$$
G=\overline{(\ZZ/2\ZZ)_F} \subset\Aut_{E/S}
$$
is called the \emph{effective model of the sign involution}. This closure can be computed locally, over the discrete valuation rings $R=\O_{S,s}$,
and we see from the discussion in Section \ref{Effective models} that $G\ra S$ is a family of group schemes of order two.
Likewise, the scheme of $G$-fixed points $\Fix(E/S)$  and the scheme of non-smoothness $\Sing(E/S)$ can   be computed locally.  
The following terminology is crucial for the whole paper:

\begin{definition}
\mylabel{admissible for kummer}
The elliptic curve $E_F$ is  called \emph{admissible for the Kummer construction} if the following three conditions hold: 
\begin{enumerate}
\item 
The intersection  $\Sing(E/S)\cap \Fix(E/S)$ is  empty.
\item
For each   $s \in S$ of bad reduction, the residue field $\kappa(s)$ has characteristic two, 
and the reduction is additive.
\item 
For each closed point $s\in S$ of characteristic two, the fiber  $\Fix(E/S)_s$ is geometrically disconnected.
\end{enumerate}
\end{definition}

Note that if $s\in S$ is a point of characteristic two and  good reduction, the fiber $\Fix(E/S)_s$
coincides with $E_s[2]$, so condition (iii) means that
the elliptic curve $E_s$ is    \emph{ordinary},  or equivalently has invariant $j\neq 0$. 

There is a variant   for pairs:
Let $E'_F$ be another elliptic curve over the field of fractions $F$, and form the Weierstra\ss{} model 
$E'$ and the effective model $G'$  of the sign involution  as above.
 
\begin{definition}
\mylabel{pairs for kummer}
The pair of elliptic curves $E_F, E'_F$ is  called \emph{admissible for the Kummer construction} if
they have the same points of bad reduction,  each elliptic curve satisfies   (i) and (ii) of the previous definition,
and furthermore the following conditions hold:
\begin{enumerate}
\item[(iii*)]
For each  point $s\in S$ of characteristic two,  at least one of the fibers $\Fix(E/S)_s$ and $\Fix(E'/S)_s$ is geometrically disconnected.
\item[(iv)]
The relative group schemes $G$ and $G'$  are isomorphic.
\end{enumerate}
\end{definition}
 
Of course, if the single elliptic curve $E_F$ is admissible for the Kummer construction, then the pair given by $E'_F=E_F$ is admissible.
Conversely,  if a pair $E_F$ and $E'_F$ is admissible  and there is at most one point $s\in S$ of characteristic two,
then at least one member of the pair is admissible.

\emph{In what follows we assume that $E_F, E'_F$ is a pair of  elliptic curves that  is admissible for the Kummer construction.}
Using condition (iv) we choose an identification  $G=G'$ between the effective models for the sign involutions
and consider the resulting diagonal $G$-action on $A=E\times_S E'$. Let us  first summarize what is known about 
the fiber-wise   quotients:

\begin{proposition}
\mylabel{fiberwise quotients}
Let  $s:\Spec(k)\ra S$ be a geometric point of characteristic $p\geq 0$.
Then the categorical quotient $A_s/G_s$ is a normal K3 surface,
and its configuration of rational double points is given by the following table:
$$
\begin{array}{llll}
\toprule
\text{\rm characteristic}	& \text{\rm effective model $G_s$}	& \text{\rm configuration of RDP}\\
\midrule
p\neq 2			& \mu_2=\ZZ/2\ZZ		& 16 A_1\\
\midrule
p=2			& \mu_2			& 16A_1 + D_4\\
			& \alpha_2			& \text{\rm $4D_4+D_4$  \; or\; $2D_8+D_4$}\\
			& \ZZ/2\ZZ			& \text{\rm $4D_4$\hspace{29pt}  \; or\; $2D_8$}\\

\bottomrule 
\end{array}
$$
\end{proposition}

\proof
The situation for  $p\neq 2$ is classical. 
The case $p=2$ follows  from \cite{Katsura 1978}, Theorem C, compare also \cite{Shioda 1974}, Remark in \S6,
and \cite{Kondo; Schroeer 2021}, Proposition 3.2, see also \cite{Schroeer 2007}.
\qed

\medskip
Of course, one may  already form   categorical quotients $A_s/G_s$ over each schematic point $s\in S$.
By the above, these are  normal K3 surfaces  over the residue field $k=\kappa(s)$, but in case of imperfectness
the minimal resolution may not be smooth.

Next, we    form    \emph{categorical quotients} over the Dedekind scheme $S$, and obtain a commutative diagram
\begin{equation}
\label{schemes and geometric quotients}
\begin{CD}
E	@<\pr_1<<	A	@>\pr_2>> 	E'\\
@VVV		@VVqV		@VVV\\
E/G	@<<<	A/G	@>>>	E'/G.
\end{CD}
\end{equation} 
Such quotients indeed exist  as algebraic spaces, for example by \cite{Rydh 2013}, Corollary 5.4.
For each $s\in S$, the universal properties of quotients give  comparison maps
\begin{equation}
\label{comparision maps}
E_s/G_s\ra (E/G)_s\quadand A_s/G_s\ra (A/G)_s\quadand E'_s/G'_s\ra (E'/G')_s.
\end{equation} 

\begin{proposition}
\mylabel{catgorical quotients}
In the diagram \eqref{schemes and geometric quotients}, all algebraic spaces   are projective and flat over $S$.
Moreover,  the comparison maps  in  \eqref{comparision maps} are isomorphisms, for all $s\in S$.
\end{proposition}

\proof
The zero sections for the Weierstra\ss{} models   yield relatively ample invertible sheaves on $E$ and $E'$, hence $A=E\times E'$ is projective.
According to \cite{Rydh 2013}, Proposition 4.7 the categorical quotient $A/G$ must be projective as well.
On the  complement $A^0=A\smallsetminus\Fix(A/S)$ of the scheme of fixed points,
the $G$-action is free. The ensuing categorical quotient $A^0/G$ is   actually a quotient of fppf sheaves, and
commutes with arbitrary base-changes, see for example \cite{Rydh 2013}, Theorem 2.16. 
It  comes with an open embedding $A^0/G\subset A/G$, whose complement
is the image $q(\Fix(A/G))$.
We see that the canonical maps $A_s/G_s\ra (A/G)_s=V_s$ are finite and birational, and are isomorphisms
outside a finite set of closed points. Using that $V_s$ has no embedded components, one gets a short exact sequence
$$
0\lra \O_{V_s}\lra \O_{A_s/G_s}\lra \shF_s\lra 0
$$
for some coherent sheaf $\shF_s$ with finite support and Euler characteristic 
$\chi(\shF_s)=\chi(\O_{A_s/G_s}) - \chi(\O_{V_s})$. We have $ \chi(\O_{A_s/G_s})=2$ because the fiber-wise quotients
are normal K3 surfaces. The function     $s\mapsto \chi(\O_{V_s})$ is constant by flatness of $V$. The generic
value is $\chi(\O_{V_\eta})=2$, because the formation of rings of invariants commutes with flat base-change.
This shows $\chi(\shF_s)=0$. In turn, the sheaves $\shF_s$ must vanish and the comparison maps  $A_s/G_s\ra V_s$ are isomorphisms.
The arguments for the elliptic curves $E$ and $E'$ are likewise.
\qed

\medskip
Note that the above reasoning on base-change of quotients appears to  be useful in other contexts as well. 
Combining the preceding two propositions, we arrive at:

\begin{theorem}
\mylabel{quotient is family of k3}
The categorical quotient $V=A/G$ is a family of normal K3 surfaces over the Dedekind scheme $S$,
and the structure morphism $V\ra S$ is projective.
\end{theorem}
 
 
In the proof for Proposition \ref{fiberwise quotients} we already described the singularities in the geometric fibers for $V\ra S$.
This immediately gives:

\begin{proposition}
\mylabel{images fixed and singular}
As a closed set, the image of the canonical map
$$
\Fix(A/S)\cup\Sing(A/S)\lra A/G=V
$$
coincides with the   locus of non-smoothness $\Sing(V/S)$.
\end{proposition}

Note that $\Fix(A/S)=\Fix(E/S)\times\Fix(E'/S)$ is finite and flat over $S$,
whereas the closed set  
$$
\Sing(A/S)=\Sing(E/S)\times\Sing(E')\subset A
$$
comprises only finitely many points $a_1,\ldots,a_r$. Their images   
$s_1,\ldots,s_r\in S$ are precisely the closed points of characteristic two where the $G$ has an infinitesimal fiber.
We like to call the images  $v_1,\ldots,v_r\in V$ the \emph{critical points}.
The residue field extensions $\kappa(s_i)\subset\kappa(v_i)$ are purely inseparable. In case of equality,
and the corresponding singularity on $V\otimes\kappa(s_i)^\perf$ is a rational double point of type $D_4$ or $B_3$,
according to \cite{Schroeer 2007}, Proposition 5.3.
In loc.\ cit., they are called \emph{singularities coming from the quadruple point}, since they stem 
from the rational point of multiplicity four on $(E\times E')\otimes\kappa(s_i)^\perf$.
For imperfect residue fields $\kappa(s_i)$, however, the  local rings
$\O_{V,v_i}$   may also be of type $G_2$ or  $A_1$,   or  be  regular.  

In contrast, the scheme $\Fix(A/S)$   yields families of singularities inside the family of normal K3 surfaces $V=A/G$.
It turns out that they are easier to understand.
Following \cite{Schroeer 2007}, Section 5 we refer to  the points $v\in V$ lying in the image of $\Fix(A/S)$ as the
\emph{singularities stemming from the fixed points}. This is indeed justified:

\begin{proposition}
\mylabel{singularities from fixed points}
For every $v\in V$ belonging to the image of $\Fix(A/S)$, the   local ring $\O_{V,v}$ is singular.
\end{proposition}

\proof
It suffices to check that the two-dimensional complete local ring $R=\O^\wedge_{V_s,v}$ is singular, where $s\in S$ is the image of $v\in V$.
Set $k=\kappa(s)$ and $k'=\kappa(v)$.
By construction, $R$ is the ring of invariants for some   action of the group $\ZZ/2\ZZ$ on a formal power series ring $A=k'[[u_1,u_2]]$.
The action is linear over $k$, and free outside the closed point of $\Spec(A)$.  As in \cite{Lorenzini; Schroeer 2020}, Proposition 3.2,
one sees that the ring of invariants is not regular. 
\qed

\section{Blowing-up  centers with embedded components}
\mylabel{Blowing ups}
 
We keep the set-up as in the preceding section, and form the family of normal K3 surfaces $V=(E\times E')/G$
over the Dedekind scheme $S$,
for a pair of elliptic curves $E_F$, $E'_F$  that is admissible for the Kummer construction,
and the resulting effective model $G$ for the sign involution.
\emph{Here we also assume that the function field $F=\O_{S,\eta}$ has characteristic zero.} 
Let $s_1,\ldots,s_r\in S$ be the finitely many points where the effective model $G$ has 
unipotent fiber, and set $F_i=\Frac(\O_{S,s_i}^\wedge)$. 
The goal of this section is to establish the following:

\begin{theorem}
\mylabel{simultaneous partial resolution}
Suppose the finite \'etale  group schemes $E_F[2]$ and $E'_F[2]$ of degree four become constant over all $F_i$, for $1\leq i\leq r$.
Then there is a simultaneous partial resolution $r:Y\ra V$ such that for each geometric point $s:\Spec(k)\ra S$
the induced map $Y_s\ra V_s$ coincides with the minimal resolution of singularities for the singularities
stemming from the fixed points.
\end{theorem}

In other words, $\Sing(Y/S)$ comprises only the  critical points $y_1,\ldots,y_r\in Y$, which correspond to
the finitely many closed points $s_1,\ldots,s_r\in S$ where $G$ becomes infinitesimal.
We shall see that outside  these $s_i$, the family $Y$ can be obtained by a blowing-up,
where the center is a family of one-dimensional schemes that are not Cartier, a technique already used in \cite{Schroeer 2007}, Sections 10 and 11. 
Here this will   involve curves with embedded components, a novel feature which might be useful in other contexts as well.
The proof requires some preparation and will be completed at the end of this section.

The projection $\pr_1:E\times E'\ra E$ induces a fibration 
$V\ra E/G$,  which can be seen as  a family of genus-one fibrations parametrized by $S$, and will be denoted by the same symbol.  Write
$\widetilde{\Fix}(E/S)\subset E/G$
for the schematic image of the fixed scheme $\Fix(E/S)\subset E$, and form  the fiber product
$$
\begin{CD}
\tilde{Z}		@>>>	V\\
@VVV			@VV\pr_1 V\\
\widetilde{\Fix}(E/S)	@>>>	E/G.
\end{CD}
$$
We now use the reduced closed subscheme  $Z=(\tilde{Z})_\red$  on the family of normal K3 surfaces  
$V=(E\times E')/G$ as center for some blowing-up 
\begin{equation}
\label{first blowing-up curve family}
\Bl_Z(V)\lra V.
\end{equation} 
The projection $Z\ra S$ is flat of relative dimension one, and   the total space is reduced and thus satisfies Serre Condition $(S_1)$. 
Note that the closed fibers $Z_s$, however,  may acquire embedded components.  

To see that this behavior is very common, just identify two   closed points in  the relative projective line $\PP^1_{\CC[[t]]}$.
The resulting closed fiber has as reduction a nodal genus-one curve,  and the sheaf of nilradicals is given by the residue
field of the singularity. The following crucial observation takes advantage of this phenomenon:

\begin{lemma}
\mylabel{center with embedded components}
Let $s:\Spec(k)\ra S$ be a geometric point such that   $V_s$ contains no $D_8$-singularity.
Then $\Bl_Z(V)_s\ra V_s$ coincides with the blowing-up of the singularities stemming from the fixed points.
\end{lemma}

\proof
It suffices to treat the case that $S$ is the spectrum of a complete discrete valuation ring $R$
with algebraically closed residue field $k=R/\maxid$, and that $s=\sigma$ is the closed point.
Replacing $R$ by a finite extension, we may also assume
that the finite \'etale scheme $E_F[2]$ of degree four over $F=\Frac(R)$ is constant.
Thus 
$\Fix(E/S)$ has four irreducible components $\Fix_i(E/S)$, $1\leq i\leq 4$. Write 
$$
\widetilde{\Fix}_i(E/S)\subset\widetilde{\Fix}(E/S)\quadand
\tilde{Z}_i\subset\tilde{Z}\quadand
Z_i\subset Z
$$
for the ensuing irreducible components.
Then $\widetilde{\Fix}_i(E/S)$ is a section for $V\ra E/G=\PP^1_R$, and     $\tilde{Z}_i\subset V$ is the corresponding
family of fibers. These fibers are irreducible, but come with multiplicity two.
In fact, $\tilde{Z}_i=\PP^1_R\oplus\O_{\PP^1_R}(-2)$
are   \emph{ribbons} in the sense of Bayer and Eisenbud \cite{Bayer; Eisenbud 1995}.
Each $Z_i\subset V$ is the ensuing family of half-fibers, and thus $Z_i=\PP^1_R$.
 
Suppose first that the group scheme $G_\sigma$ is multiplicative.   Then the  $\Fix_i(E/S)$    are pairwise disjoint, and likewise
$Z_i=\PP^1_R$ are pairwise disjoint families of half-fibers with respect to $V\ra E/G=\PP^1_R$.
Using \cite{Schroeer 2007}, Proposition 10.5 we infer that $\Bl_Z(V)\ra V$ gives fiberwise the minimal resolution of singularities,  
and the assertion follows.

It remains to treat the case that $G_\sigma$ is unipotent. Since $V_\sigma$ contains no $D_8$-singularity,
 $\widetilde{\Fix}(E/S)_\sigma=\{a,b\}$ comprises two $k$-rational points,  each with multiplicity two.
Without loss of generality $\widetilde{\Fix}_1(E/S)\cap \widetilde{\Fix}_2(E/S)=\{a\}$. In turn, $(Z_1\cup Z_2)_\sigma$ contains the fiber 
$V_a=V\times_{E/G}\Spec \kappa(a)=\PP^1_k\oplus\O_{\PP^1_k}(-2)$,
and the inclusion $V_a\subset (Z_1\cup Z_2)_\sigma$ is an equality outside $\Sing(V/S)$. From
$$
\chi(\O_{(Z_1\cup Z_2)_\sigma}) =  \chi(\O_{(Z_1\cup Z_2)_\eta}) = 2\quadand \chi(V_a)=\chi(\O_{\PP^1_k})\oplus \O_{\PP^1_k}(-2))=0
$$
we see that in the short exact sequence $0\ra \shI\ra \O_{(Z_1\cup Z_2)_\sigma}\ra \O_{V_a}\ra 0$
the kernel is a skyscraper sheaf with $h^0(\shI)=2$, which is supported by the two $D_4$-singularities $z_1,z_2\in V_\sigma$  lying on the fiber $V_a=\pr_1^{-1}(a)$.
Since  $V_a\subset V$ is Cartier, and $Z_1\cup Z_2\subset V$ is not Cartier at each of the  $A_1$-singularities on $V_\eta$ specializing to $z_1$ or $z_2$,
it follows from the Nakayama Lemma 
that $\shI$ is non-zero at both $z_1,z_2$. In light of $h^0(\shI)=2$ we must have  $\shI_{z_i}\simeq\kappa(z_i)$.
In turn, the blowing-up $Y\ra \Bl_Z(V)$ becomes near $z_i\in V_\sigma$ 
the blowing-up of an ideal of the form $\maxid\cdot\ideala\subset\O_{V_\sigma,z_i}$ for some invertible   $\ideala$.
This result is the same  as the  blowing-up of $\maxid\subset\O_{V_\sigma,z_i}$. 
Summing up, $Y_\sigma\ra V_\sigma$ coincides with the blowing-up of the singularities $z\in V_\sigma$, endowed with reduced scheme structure.
\qed

\medskip
We continue to examine the above blowing-up  $\Bl_Z(V)$, \emph{now under the additional assumption
that  the finite \'etale group schemes $E_\eta[2]$ and $E'_\eta[2]$ of degree four are constant}.
Then all singularities in $V_\eta$ are $F$-valued, and we write $\Delta_1,\ldots,\Delta_{16}$
for the exceptional divisors on the resolution of singularities $\Bl_Z(V)_\eta$, arranged
in an arbitrary order. Write
$$
\Bl_{Z,16}(V)=V_{16}\lra V_{15}\lra \ldots\lra V_1\lra V_0=\Bl_Z(V)
$$
for the sequence of blowing-ups, where for $V_i=\Bl_{Z_i}(V_{i-1})$  
the center is taken as the schematic closure of the $\Delta_i\subset (V_{i-1})_\eta=\Bl_Z(V)_\eta$.
Of course, such a center may already be Cartier, in which case the blowing-up becomes an identity.

\begin{proposition}
\mylabel{sixteen blowing-ups}
For each geometric point $s:\Spec(k)\ra S$ such that   $V_s$ contains no $D_8$-singularity, the composite morphism
$\Bl_{Z,16}(V)_s\ra V_s$ is the   minimal resolution of all singularities stemming from the fixed points, and an isomorphism otherwise.
\end{proposition}
 
\proof
It suffices to treat the case that $S$ is the spectrum of a complete discrete valuation ring, with closed point
$\sigma=s$. In light of Proposition \ref{center with embedded components}, there is nothing to prove if $G_\sigma$
is multiplicative. So we assume that the residue field $k$ has characteristic two, and that $G_\sigma$ is unipotent.

Write $r_s:X_s\ra V_s$ for the minimal resolutions of singularities.
The exceptional divisors on   $X_\eta$, together with the strict transform of the half-fibers over $(E/G)_\eta=\PP^1_\eta$,
take the form:
\begin{equation}
\label{dual graphs for i0^* fibers}
\begin{gathered}
\begin{tikzpicture}
[node distance=1cm, font=\small]
\tikzstyle{vertex}=[circle, draw, fill=white, inner sep=0mm, minimum size=1.5ex]

\node[vertex,fill=black]	(C2)  	at (1,0) [  label=above:{$\Phi_1$}]{}; 
\node[vertex]	(C0)	[above left of=C2, label=left:{$\Upsilon_{1,1}$}]{};
\node[vertex]	(C1)	[below left of=C2, label=left:{$\Upsilon_{1,2}$}]{};
\node[vertex]	(C3)	[above right of=C2, label=right:{$\Upsilon_{1,3}$}]{};
\node[vertex]	(C4)	[below right of=C2, label=right:{$\Upsilon_{1,4}$}]{};

\draw [thick] (C3)--(C2)--(C4);
\draw [thick] (C0)--(C2)--(C1);

\node[] (E)		at (3.5,0) [label=right:{$+\quad\ldots\quad+$}]{};

\node[vertex,fill=black]	(D2)  	at (9,0) [  label=above:{$\Phi_4$}]{}; 
\node[vertex]	(D0)	[above left of=D2, label=left:{$\Upsilon_{4,1}$}]{};
\node[vertex]	(D1)	[below left of=D2, label=left:{$\Upsilon_{4,2}$}]{};
\node[vertex]	(D3)	[above right of=D2, label=right:{$\Upsilon_{4,3}$}]{};
\node[vertex]	(D4)	[below right of=D2, label=right:{$\Upsilon_{4,4}$}]{};

\draw [thick] (D3)--(D2)--(D4);
\draw [thick] (D0)--(D2)--(D1);

\end{tikzpicture}
\end{gathered}
\end{equation} 
Here the black vertex signifies the strict transform of the half-fiber, 
and the $\Upsilon_{i,j}$ correspond to the $\Delta_1,\ldots,\Delta_{16}$.
The exceptional divisors on   $X_\sigma$ form two configurations,    with  the following dual graphs and $i=1,2$:
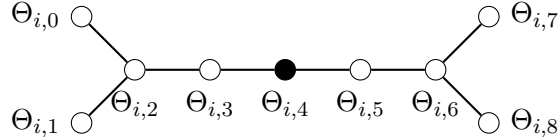
\begin{equation*}
\label{I12*}
\begin{gathered}
\begin{tikzpicture}
[node distance=1cm, font=\small]
\tikzstyle{vertex}=[circle, draw, fill=white, inner sep=0mm, minimum size=1.5ex]
\node[vertex]	(C2)  	at (1,0) 	[label=below:{$\Theta_{i,2}$}] 		{};
\node[vertex]	(C3)			[right of=C2, label=below:{$\Theta_{i,3}$}]	{};
\node[vertex,fill=black]	(C4)			[right of=C3, label=below:{$\Theta_{i,4}$}]	{};
\node[vertex]	(C5)			[right of=C4, label=below:{$\Theta_{i,5}$}]	{};
\node[vertex]	(C6)			[right of=C5, label=below:{$\Theta_{i,6}$}]	{};
 
\node[vertex]	(C7)			[above right of=C6, label=right:{$\Theta_{i,7}$}]	{};
\node[vertex]	(C8)			[below right of=C6, label=right:{$\Theta_{i,8}$}]	{};
\node[vertex]	(C0)			[above left  of=C2, label=left :{$\Theta_{i,0}$}]	{};
\node[vertex]	(C1)			[below left  of=C2, label=left :{$\Theta_{i,1}$}]	{};
 
\draw [thick] (C2)--(C3)--(C4)--(C5)--(C6);
\draw [thick] (C0)--(C2)--(C1);
\draw [thick] (C8)--(C6)--(C7);
\end{tikzpicture}
\end{gathered}
\end{equation*}
Recall that $\Bl_Z(V)\ra V$ is the blowing-up of the schematic image of $\Phi_1+\ldots+\Phi_4$.
The  degenerate fibers in $\Bl_Z(V)_\sigma\ra (E/G)_\sigma=\PP^1_k$ are the images of
$$
\Theta_{1,2}+\Theta_{1,4}+\Theta_{1,6}\quadand \Theta_{2,2}+\Theta_{2,4}+\Theta_{2,6}.
$$
The remaining singularities on $\Bl_Z(V)_\sigma$ are merely $A_1$-singularities. 
We also see that each $\Upsilon_{i,j}\subset X_\eta=(\Res_{V/S})_\eta$
specializes to the schematic image of some of the remaining $\Theta_{r,s}$ with 
$r=1,2$ and $s=0,1,3,5,7,8$. In turn, the blowing-up of the schematic image of  $\Upsilon_{i,j}$ resolves the
singularities  that lie on $\Theta_{r,s}$. The assertion follows, because the $\Upsilon_{ij}$ are exactly
the $\Delta_1,\ldots,\Delta_{16}$.
\qed

\medskip
\emph{Proof of Theorem \ref{simultaneous partial resolution}.}
We proceed by induction on the number $r\geq 0$ of points $\sigma_1,\ldots,\sigma_r\in S$ where
the effective model $G$ is unipotent. For $r=0$, the group scheme $G$ is multiplicative,
the singularities $v\in V_s$ in the geometric fibers stemming from the fixed points are 
of type $A_1$, and the assertion follows from Lemma \ref{center with embedded components}.

Suppose now $r\geq 1$, and that the assertion holds for $r-1$.
Let $U\subset S$ be the complement of $\sigma=\sigma_r$.
By the induction hypothesis, we find the desired simultaneous partial resolution $r_U:Y_U\ra V_U$.
Suppose we also have the desired simultaneous partial resolution over   $\O^\wedge_{S,\sigma}$.
According to \cite{Artin 1970}, Theorem 3.2 this formal modification already exists over $\O_{S,\sigma}$, and thus 
extends to some $r_{U'}:Y_{U'}\ra V_{U'}$ for some  affine open neighborhood $U'$ of $\sigma\in S$.
Shrinking if necessary, we may assume that it coincides with $Y_U$ on the overlap $U\cap U'$.
The desired $r:Y\ra V$ is then obtained by gluing.

This reduces us to the case that $S$ is the spectrum of a complete discrete valuation ring $R$,
with closed point $\sigma\in S$, residue field $k=\kappa(\sigma)$, and unipotent $G_\sigma$.
If the geometric closed fiber contains no $D_8$-singularities, the assertion follows from Proposition \ref{sixteen blowing-ups}.
It remains to treat the case that there is a $D_8$-singularity. Recall that by assumption, the finite \'etale group schemes
$E_F[2]$ and $E'_F[2]$ are constant.
Without loss of generality we may assume that $\Fix(E/S)$ is connected and $\Fix(E'/S)$ is   disconnected.

Write $X_s\ra V_s$, $s\in S$ for the minimal resolution of singularities.
The exceptional divisors on $X_\sigma$, together with the strict transform $\Theta_8$ of the corresponding half-fiber
over $(E/G)_\sigma=\PP^1_\sigma$,
form a configuration with the following dual graph:
\begin{equation*}
\label{I11*}
\begin{gathered}
\begin{tikzpicture}
[node distance=1cm, font=\small]

\tikzstyle{vertex}=[circle, draw, fill=white, inner sep=0mm, minimum size=1.5ex]
\node[vertex]	(C2)  	at (1,0) 	[label=below:{$\Theta_2$}] 		{};
\node[vertex]	(C3)			[right of=C2, label=below:{$\Theta_3$}]	{};
\node[vertex]	(C4)			[right of=C3, label=below:{$\Theta_4$}]	{};
\node[vertex]	(C5)			[right of=C4, label=below:{$\Theta_5$}]	{};
\node[vertex]	(C6)			[right of=C5, label=below:{$\Theta_6$}]	{};
\node[vertex]	(C7)			[right of=C6, label=below:{$\Theta_7$}]	{};
\node[vertex,fill=black]	(C8)			[right of=C7, label=below:{$\Theta_8$}]	{};
\node[vertex]	(C9)			[right of=C8, label=below:{$\Theta_9$}]	{};
\node[vertex]	(C10)			[right of=C9, label=below:{$\Theta_{10}$}]	{};
\node[vertex]	(C11)			[right of=C10, label=below:{$\Theta_{11}$}]	{};
\node[vertex]	(C12)			[right of=C11, label=below:{$\Theta_{12}$}]	{};
\node[vertex]	(C13)			[right of=C12, label=below:{$\Theta_{13}$}]	{};
\node[vertex]	(C14)			[right of=C13, label=below:{$\Theta_{14}$}]	{};

\node[vertex]	(C15)			[above right of=C14, label=right:{$\Theta_{15}$}]	{};
\node[vertex]	(C16)			[below right of=C14, label=right:{$\Theta_{16}$}]	{};
\node[vertex]	(C0)			[above left  of=C2, label=left :{$\Theta_0$}]	{};
\node[vertex]	(C1)			[below left  of=C2, label=left :{$\Theta_1$}]	{};
 
\draw [thick] (C2)--(C3)--(C4)--(C5)--(C6)--(C7)--(C8)--(C9)--(C10)--(C11)--(C12)--(C13)--(C14);
\draw [thick] (C0)--(C2)--(C1);
\draw [thick] (C15)--(C14)--(C16);

\end{tikzpicture}
\end{gathered}
\end{equation*}
As before,  the black vertex signifies the strict transform of the half-fiber.
The exceptional divisors on $X_\eta$, together with the strict transform of the half-fibers over $(E/G)_\eta=\PP^1_\eta$
are already given in \eqref{dual graphs for i0^* fibers}.

We now construct iterated blowing-ups of $Y_4\ra \ldots \ra Y_1\ra V$, and it is most convenient to specify the
centers as schematic images of certain combinations of $\Phi_i$ and $\Upsilon_{i,j}$  with respect to the canonical morphism $X_\eta\ra Y_i$.
We   proceed in four steps:

\medskip
\textbf{Step 1.} Let $Y_1=\Bl_Z(V)$ be the blowing-up   where the center $Z\subset V$
is the schematic image of $\Phi_1\cup\Phi_2$. This pair of curves specializes to the schematic images
of $2\Theta_8$, with embedded components at the singularities.
The situation  is essentially the same  as in Lemma \ref{center with embedded components}:  The degenerate fiber in
$Y_{1,\sigma}\ra (E/G)_\sigma=\PP^1_\sigma$ comprises the images of 
$$
\Theta_6+\Theta_8+\Theta_{10},
$$
and the singularities are $(D_6+A_1)+(A_1+D_6)$,  contained in the images of $\Theta_6+\Theta_{10}$.
(For the behavior of rational double points under standard blowing-ups, compare
\cite{Schroeer 2009}, proof of Proposition 4.2).

\medskip
\textbf{Step 2.} 
Let $Y_2=\Bl_{Z_1}(Y_1)$ be the blowing-up where the center $Z_1\subset Y_1$ is the schematic image
of $\Upsilon_{1,1}+\ldots +\Upsilon_{1,4}$. These four curves specialize in pairs to the schematic images of $2\Theta_6$ and $2\Theta_{10}$,
with embedded components at the four singularities.
Again the situation is analogous to Lemma \ref{center with embedded components}, and  the degenerate fiber in 
$Y_{2,\sigma}\ra(E/G)_\sigma=\PP^1_\sigma$ comprises the image of 
$$
\Theta_4+ \left(\sum_{i=6}^{10}\Theta_i\right)+\Theta_{12}.
$$
Now the singularities are $(D_4+A_1)+(A_1+D_4)$, which are contained in the images of $\Theta_4+\Theta_{12}$.

\medskip
\textbf{Step 3.} 
Let $Y_3=\Bl_{Z_2}(Y_2)$ be the blowing-up where the center $Z_2\subset Y_2$ is the schematic image
of $\Upsilon_{2,1}+\ldots +\Upsilon_{2,4}$. These four curves specialize pairwise to the images of $2(\Theta_4+\Theta_6+\Theta_7)$ and 
$2(\Theta_9+\Theta_{10}+\Theta_{12})$, with embedded components at the four singularities.
Again the situation is analogous to Lemma \ref{center with embedded components}. The degenerate fiber in $Y_{3,\sigma}\ra\PP^1_\sigma$
comprises the images of
$$
\Theta_2+\left(\sum_{i=4}^{12}\Theta_i\right) + \Theta_{14},
$$
and the singularities are $(A_1+A_1+A_1)+ (A_1+A_1+A_1)$, contained in the images of $\Theta_2+\Theta_{14}$.
After changing the enumeration of the $\Upsilon_{2,i}$, we may assume
that the image of  $\Upsilon_{2,1}$ specializes to the image of $\Theta_2$,
and the image of $\Upsilon_{2,4}$ specializes to the image of $\Theta_{14}$,
now without any embedded components.

\medskip
\textbf{Step 4.} Let $Y=\Bl_{Z_3}(Y_3)$ be the blowing-up where the center $Z_3\subset Y_3$ is the schematic image
of $\Upsilon_{2,1} +\Upsilon_{2,4}$.  It then follows from \cite{Schroeer 2007}, Proposition 10.5  that  the composite map $Y\ra V$
is indeed the desired simultaneous partial resolution of singularities.
\qed

\section{Monodromy representations}
\mylabel{Monodromy}

Let $S$ be the spectrum of a    discrete valuation ring $R$ whose field of fractions $F=\Frac(R)$ has
characteristic zero, and whose residue field $k=R/\maxid$ has characteristic two.
Write $\eta,\sigma\in S$ for the generic and closed point, respectively.
Let $E_\eta$, $E'_\eta$ be a pair of elliptic curves over $F=\kappa(\eta)$ that are admissible for the Kummer construction
and form the   family of normal K3 surfaces $V=(E\times E')/G$, where $G$ is the effective model of the sign involution.
Write 
$$
r_\eta:X_\eta\lra V_\eta\quadand r_\sigma:X_\sigma\lra V_\sigma
$$
for the fiber-wise minimal resolutions of singularities. Then $X_\eta$ is a K3 surface;
the same is true for $X_\sigma$ provided that $k$ is perfect.
In this section we introduce various $\ell$-adic sheaves on the spectra of $F$ and $k$ that
encode crucial information, 
and establish both necessary and sufficient conditions for the existence of a simultaneous minimal resolution.
 
\emph{Throughout, we assume that there is a simultaneous partial resolution $Y\ra V$
that fiberwise resolves all singularities stemming from the fixed points}. 
Note that by Lemma \ref{center with embedded components}, this   automatically holds  if the effective model $G$ is multiplicative.
The assumption is  mainly for the sake of exposition, and could be avoided at the expense of using $\ell$-adic cohomology with compact supports.

Write $g:Y\ra S$ and $f_s:X_s\ra\Spec \kappa(s)$, $s\in S$ for the structure morphisms, which sit in  commutative diagrams
$$ 
\begin{tikzcd}
X_s\ar[rr,"r_s"]\ar[dr,"f_s"']	&	& Y_s\ar[dl,"g_s"]\\	 
			& S.
\end{tikzcd}
$$
As customary, set $A=E\times E'$  and let $q:A\ra A/G=V$ be the quotient map.
We then get a   diagram 
\begin{equation}
\label{various maps for kummer}
\begin{CD}
\uH^1(E_\eta,\ZZ_\ell)\otimes \uH^1(E'_\eta,\ZZ_\ell)	@>\cup >>	\uH^2(A_\eta,\ZZ_\ell)\\
@.						@AAq_\eta^* A\\
@.						\uH^2(V_\eta,\ZZ_\ell)	@>>r_\eta^*> H^2(X_\eta,\ZZ_\ell)
\end{CD}
\end{equation} 
of $\ell$-adic systems on the scheme $\Spec(F)$. Here $\ell>0$ is an odd prime, and the map on the left stems from the cup product via
$\alpha\otimes \alpha'\mapsto \pr_1^*(\alpha)\cup\pr_2^*(\alpha')$.

\begin{proposition}
\mylabel{induced maps on cohomology}
In the above diagram, both horizontal maps are injective, and the vertical map is bijective.
\end{proposition}

\proof
The cup product is injective by the K\"unneth Theorem. The first higher direct images $R^1r_{\eta,*}(\ZZ_\ell)$ vanishes,
because the exceptional divisors are projective lines. Now fix an algebraic closure $F^\alg$.
The Leray--Serre spectral sequence gives an exact sequence
\begin{equation}
\label{leray serre for resolution}
0\lra H^2(V\otimes F^\alg,\ZZ_\ell)\stackrel{q^*}{\lra} H^2(X\otimes F^\alg,\ZZ_\ell)\lra H^0(V_\eta,R^2r_{\eta,*}(\ZZ_\ell)).
\end{equation} 
This already yields the injectivity of $r_\eta^*$. We also see that  the arrow on the right has finite cokernel, because
the intersection form on   exceptional divisors is negative definite. Consequently
$H^2(V\otimes F^\alg,\ZZ_\ell)$ has rank $22-16=6$. Moreover, its cup product is non-degenerate,
in light of Poincar\'e Duality for $X\otimes F^\alg$,  with discriminant   a 2-power, because the exceptional curves
are $(-2)$-curves.
Bijectivity of $q^*_\eta$ is now clear, because the groups have the same rank, and the discriminant is an $\ell$-adic unit.
\qed

\medskip
We now introduce certain $\ell$-adic subsystems
$$
\Psi_\eta \subset \uH^2(X_\eta,\ZZ_\ell)\quadand \Psi_\sigma\subset \uH^2(X_\sigma,\ZZ_\ell).
$$
The former is the  image of the composition stemming from  diagram \eqref{various maps for kummer},  formed with the inverse of $q_\eta^*$.
The latter is defined   in terms of its Tate twist, by letting
$\Psi_\sigma(1)\subset \uH^2(X_\sigma,\ZZ_\ell)(1)= \uH^2(X_\sigma,\ZZ_\ell(1))$
be the subsheaf generated by the cycle classes that are supported by the exceptional locus $\Exc_{X_\sigma/V_\sigma}$
defined over  finite separable extensions of $\kappa(\sigma)$. Obviously
$\rank(\Psi_\eta)=2\cdot 2=4$ and $\rank(\Psi_\sigma)\leq 4$.

\begin{lemma}
\mylabel{local systems match}
Suppose there is a simultaneous resolution of singularities $r:X\ra V$, for some family $f:X\ra S$ of K3 surfaces.
Then there is a local subsystem 
$$
\Psi_S\subset f_*(\ZZ_{\ell,X})
$$
whose restrictions to   $s\in S$ coincide with    $\Psi_s$ introduced above, and this has rank four.
\end{lemma}
 
\proof
One may view the Tate twist $\Psi_\eta(1)\subset \uH^2(X_\eta,\ZZ_\ell(1))$ as the orthogonal complement
of the sheaf of cycle classes of those  curves that are vertical with respect to the two projections 
of $X_\eta=E_\eta\times E'_\eta$. Likewise, 
$\Psi_\sigma(1)$ is the orthogonal complement stemming from the two projections of $X_\sigma=E_\sigma\times E'_\sigma$,
but now only with the cycle classes not stemming from the critical point.
In both cases, the number of cycle classes is $18=16+2$, and all belong to irreducible curves mapping to
the image of 
$$
\Fix(E/S)\ra E/G\quad\text{or}\quad \Fix(E'/S)\ra E'/G.
$$
Since the fibrations  and the fixed loci exist  over the base, so does the orthogonal complement  $\Psi_S$.
Then  $\rank(\Psi_s)=\rank(\Psi_\eta)=4$.
\qed

\medskip
We now translate the situation into representation theory and eventually into linear algebra.
Fix separable  closures $k^\sep$ and $F^\sep$. The former yields the strict henselization $R^\sh$.
Then $F^\sh=\Frac(R^\sh)$ is a separable extension of $F$, and we choose an $F$-embedding $F^\sh\subset F^\sep$.
With respect to these field extensions we  now form the Galois groups $\Gamma_\eta=\Gal(F^\sep/F)$ and $\Gamma_\sigma=\Gal(k^\sep/k)$,
and consider  the   monodromy representations
$$
\rho_\eta:\Gamma_\eta\lra \GL (\Psi_\eta(F^\sep))\quadand  \rho_\sigma:\Gamma_\sigma\lra \GL(\Psi_\sigma(k^\sep)).
$$
By definition, the former is a \emph{Kronecker product} $\rho_\eta=\epsilon_\eta\otimes\epsilon'_\eta$,
where the tensor factors are the monodromy representations
\begin{equation}
\label{monodromy representations elliptic curves}
\epsilon_\eta:\Gamma_\eta\lra \GL( H^1(E\otimes F^\sep,\ZZ_\ell))\quadand \epsilon'_\eta:\Gamma_\eta\lra \GL(H^1(E'\otimes F^\sep,\ZZ_\ell))
\end{equation} 
stemming from the elliptic curves $E_\eta$ and $E'_\eta$.
Finally, write  
$$
\chi_\sigma:\Gamma_\sigma\lra\Aut \left(\bigcup_{\nu\geq0}\mu_{\ell^\nu}(k)\right) = \invlim_{\nu\geq 0}(\ZZ/\ell^\nu\ZZ)^\times=\ZZ_\ell^\times
$$
for the $\ell$-adic \emph{cyclotomic character}. Its kernel corresponds to the subfield $k^\cyc\subset k^\sep$
generated by the $\ell^\nu$-th roots of unity, $\nu\geq 0$.
Note that the  corresponding field $F^\cyc\subset F^\sep$ is already contained in $F^\sh$.

Suppose that the minimal resolution  $X_\sigma$ is smooth, and let   $\Theta_0,\ldots,\Theta_3\subset X_\sigma\otimes k^\sep$
be the four exceptional divisors mapping to $v_\crit\in V$, indexed as in the Bourbaki tables 
(\cite{LIE 4-6}, page 256).
The cycle classes $\cl(\Theta_i)\in H^2(X_\sigma\otimes k^\sep,\ZZ_\ell(1))$ form a
$\ZZ_\ell$-basis in  $\Psi_\sigma(1)(k^\sep)$, and the monodromy  action permutes the basis members.
In turn, the Tate twist becomes a \emph{permutation representation} 
$\rho_\sigma\otimes\chi_\sigma:\Gamma_\sigma\ra S_4\subset\GL_4(\ZZ_\ell)$.
 
\begin{proposition}
\mylabel{monodromy represenations}
Suppose  there is a simultaneous minimal resolution of singularities $r:X\ra V$. Then the following holds:
\begin{enumerate}
\item 
Both $F=\Frac(R)$ and $k=R/\maxid$  contain   a primitive third root of unity. 
\item 
All three   monodromy representations $\epsilon_\eta$, $\epsilon'_\eta$ and  $\rho_\sigma$ are scalar representations,
and we have $\rho_\sigma=\chi_\sigma^{-1}$.
\item  
The restrictions of  $\epsilon_\eta$ and $\epsilon'_\eta$  to $\Gal(F^\sep/F^\text{cyc})$ are isomorphic,
and their   eigenvalues  belong to $\{\pm 1\}\subset\ZZ_\ell^\times$.
\item
The base changes $E\otimes F^\sh$ and $E'\otimes F^\sh$ acquire good reduction over a common   field extension $F^\sh\subset L$
of degree at most two.
\item 
If the  elliptic curves  $E_\eta, E'_\eta$ have the same $j$-invariant,  
then    their base-changes to $F^\sh$ are already isomorphic, at least if   $j\neq 1728$ or $ F^\sh$ contains a primitive fourth root of unity.
\end{enumerate}
\end{proposition}
 
\proof
First note that restricting a finite \'etale scheme over $S$ to its fiber over the closed point $\sigma$
gives  an equivalence of Galois categories.
Our Galois groups, which are the automorphism groups for the fiber functors for $F^\sep$ and $k^\sep$,
are thus related by a surjective homomorphism $h:\Gamma_\eta\ra \Gamma_\sigma$, compare the discussion 
in \cite{Ruelling; Schroeer 2025}, Section 2. In light of Lemma \ref{local systems match}, the representations 
$\rho_\eta=\epsilon_\eta\otimes\epsilon'_\eta$ and $\rho_\sigma\circ h$
are isomorphic. 
Let $v\in V_s$ be the critical singularity.
According to \cite{Schroeer 2007}, Proposition  5.3 this is a rational double point of type $D_4$
if the field $k$ contains a primitive third root of unity, and of type $B_3$ otherwise.
In turn, for  each $g\in \Gamma_\sigma$  the ensuing   matrix $(\rho_\sigma\otimes\chi_\sigma)(g)$ is similar to 
one of the permutation matrices
\begin{equation}
\label{two   representation matrices}
\begin{pmatrix}
E_2	&\\
	& 1	& 0\\
	& 0	& 1\\
\end{pmatrix}
\quadand
\begin{pmatrix}
E_2	&\\
	& 0	& 1\\
	& 1	& 0\\
\end{pmatrix}.
\end{equation} 
Both   are diagonalizable in characteristic $\ell\neq 2$, and
we conclude that the characteristic polynomial for $\rho_\sigma(g)$ takes the form 
$\chi_{\rho_\sigma(g)}(T)=(T-\alpha)^3(T \pm \alpha)$
with eigenvalue $\alpha=\chi^{-1}_\sigma(g)$. By Lemma \ref{kronecker product} below,
the factors in the Kronecker product  $\epsilon_\eta(g)\otimes\epsilon'_\eta(g)$ are   homotheties.
Thus actually $\rho_\sigma(g)=\chi^{-1}_\sigma(g)$, which already establishes assertion (ii).
Furthermore,  $(\rho_\sigma\otimes\chi_\sigma)(g)$
is given by the first matrix in \eqref{two   representation matrices}, which yields (i).
 
To verify (iii), suppose that $g\in \Gamma_\eta$ fixes the cyclotomic field $F^\cyc$.
Then the  respective eigenvalues $\alpha,\alpha'\in \QQ_\ell$ of the    homotheties $\epsilon(g),\epsilon'(g)$
satisfy $\alpha\alpha'=1$. On the other hand,  the Galois actions respect  the Weil pairings on   $E[\ell^\nu]$ and $E'[\ell^\nu]$,
so the   $\epsilon(g), \epsilon'(g)$ are symplectic, and thus have determinant $d=1$.
For the eigenvalues, this means $\alpha^2=\alpha'^2=1$, and therefore  $\alpha,\alpha'\in\{\pm 1\}$.
Using  $\alpha\alpha'=1$ we also get $\alpha=\alpha'$, which establishes (iii).

The cyclotomic field $F^\cyc$ is contained in $F^\sh$. The common restriction of $\epsilon$ and $\epsilon'$
define a homomorphism $\Gal(F^\sep/ F^\sh)\ra \{\pm 1\}$, so the corresponding field extension $F^\sh\subset L$
has degree at most two.
By construction, the  finite \'etale group schemes $E_F[\ell^\nu]$ and $E'_F[\ell^\nu]$ become constant over $L$,
for all $\nu\geq 0$.
According to the Criterion of N\'eron--Ogg--Shafarevich (\cite{Serre; Tate 1968}, Theorem 1), 
the elliptic curves $E_L$ and $E'_L$ have good reduction, thus (iv) holds.

It remains to verify (v), and for this we may assume that $S=\Spec(R)$ is strictly henselian,
such that $F^\sh=F$.  If the Weierstra\ss{} models $E,E'$ are smooth, the closed fibers
have the same $j$-invariant, and we find an isomorphism $h:E_\sigma\ra E'_\sigma$.
Since $R$ is henselian, we can lift it to an isomorphism $h:E\ra E'$. 

Suppose now that $E_F$, $E'_F$ have bad reduction. Let $F\subset L$ be the quadratic extension
over which they acquire good reduction, and set $\Gamma=\Gal(L/F)$. By the previous paragraph, we find an identification
$E_L=E'_L$. In other words, $E'_F$ is a twisted form of $E_F$.
We thus have $E'_F=Q\wedge^\Gamma E_F$ via some torsor $Q$ with respect to $\Gamma=\Aut_{E_F/F}$.
The latter is isomorphic to $\mu_{2n}$ with  some $1\leq n\leq 3$, according to \cite{Deligne 1975}, Proposition 5.9.
Using the norm map for the separable quadratic extension $F\subset L$,
we see that the torsor is induced from  $P$ with respect to   $\mu_2$.  The latter 
corresponds to the constant group $\Gamma=\{\pm 1\}$ given by the sign involution.
Applying Lemma \ref{fact on twisting} below with $(X,X_0,e)=(E_F,E_F[\ell],e)$ and using (iii), we infer that $E'_F\simeq E_F$.
\qed

\medskip
Let us record the following consequence, which is not obvious from Definition \ref{pairs for kummer}:

\begin{corollary}
\mylabel{potentially good reduction}
Both elliptic curves $E_\eta$ and $E'_\eta$ have potentially good reduction.
\end{corollary}

We now obtain a sufficient condition for resolutions of singularities:

\begin{theorem}
\mylabel{sufficient condition for resolution}
A simultaneous minimal resolution of singularities $r:X\ra V$ exists provided the following conditions hold:
\begin{enumerate}
\item The field $F$ contains a primitive third root of unity.
\item The purely inseparable extension $k\subset \kappa(v_\crit)$ is an equality.
\item The elliptic curves $E\otimes F^\sh$ and $E'\otimes F^\sh$ are isomorphic, and acquire good reduction 
over a quadratic extension $F^\sh\subset L$. 
\end{enumerate}
\end{theorem}

\proof
Let $r_\sigma:X_\sigma\ra Y_\sigma$ be the minimal resolution of singularities. 
Assumption (ii) ensures that $X_\sigma$ is smooth, as already remarked below Proposition \ref{images fixed and singular}.
We now seek to apply Shepherd-Barron's result (\cite{Shepherd-Barron 2021}, Corollary 2.14).
For this we have to ensure that all exceptional curves in $X_\sigma$ are isomorphic to $\PP^1_\sigma$,
an  assumption made in loc.\ cit., beginning of Section 2. This indeed holds by our assumption (i),
as in \cite{Schroeer 2007}, Proposition 5.3.

Consider the monodromy representations
$$
\varphi_\eta:G_\eta\lra \invlim \GL (E[\ell^\nu](F^\sep))\quadand \varphi'_\eta:G_\eta\lra \invlim \GL (E'[\ell^\nu](F^\sep))
$$
stemming from the Tate modules. By Hensel's Lemma, the 
restrictions of the above to $\Gal(F^\sep/F^\sh)$ factor over the quotient $\Gal(L/F^\sh)$, which is cyclic of order two.
Since $2$ is an $\ell$-adic unit, for each 
$g\in \Gal(F^\sep/F^\sh)$, the images $\varphi_\eta(g)$ and $\varphi'_\eta(g)$ are diagonalizable,
and have 
$N=(\begin{smallmatrix}\alpha&0\\0& \beta\end{smallmatrix})$
as common Jordan normal form, for some $\alpha,\beta\in \mu_2(\QQ_\ell)=\{\pm 1\}$. Since $F^\sep $ contains the cyclotomic field $F^\cyc$,
the monodromy representation respects the Weil pairing, so the above representation matrix must be symplectic,
and thus has determinant $\alpha\beta=1$. Consequently $\alpha=\beta$, so the above matrix takes the form $N=\pm E_2$.

Using \cite{Raynaud 1970}, Proposition 6.2.1 we get an identification $A_F[\ell^\nu] = \uH^1(A,\ZZ/\ell^\nu\ZZ(1))$
of $\ZZ/\ell^\nu\ZZ$-local systems
and   infer that the tensor product of the
$\ell$-adic representations \eqref{monodromy representations elliptic curves} is trivial when restricted to $F^\sep$. 
By \cite{Shepherd-Barron 2021}, Corollary 2.14 the desired simultaneous minimal resolution $r:X\ra Y$ exists.
\qed
 
\medskip
Building on an observation of Serre (\cite{Serre 1989}, page 55), Zarhin analyzed the special role of roots of unity in $\ell$-adic cohomology
\cite{Zarhin 2017}, and it would be interesting to see if there are deeper connections.

In the   proof for Proposition \ref{monodromy represenations}, we have used some auxiliary facts. The first comes from linear algebra:
Let $V$, $W$ be two non-zero finite-dimensional vector spaces over some field  $K$, and set  $m=\dim(V)$ and $n=\dim(W)$.
For $f  \in\End(V)$ and $g\in\End(W)$ we form the \emph{Kronecker product}
$f\otimes g\in \End(V\otimes W)$. Recall that the characteristic polynomials are related by
$$
\chi_f(T)=\prod_{i=1}^m(T-\lambda_i),\quad \chi_g(T)=\prod_{j=1}^n(T-\mu_j)\quadand \chi_{f\otimes g}(T)=\prod_{i,j}(T-\lambda_i\mu_j),
$$
for certain $ \lambda_i,\mu_j\in K^\alg$. Also recall that endomorphisms that become diagonalizable after some field extension
are called \emph{semi-simple}.

Obviously, if $f$ and $g$ are diagonalizable or semi-simple, the same property holds for the Kronecker product.
Using Jordan normal forms, one sees that if $f\otimes g$ is   semisimple, the same holds for the factors.
On the other hand,   $f\otimes g$ can be  diagonalizable  without $f$ and  $g$ being so:
This already happens when both $f$ and $g$ have rational normal form $(\begin{smallmatrix} 0&-1\\1&0\end{smallmatrix})$,
over the field $K=\RR$. Suitable assumptions preclude such phenomena:

\begin{lemma}
\mylabel{kronecker product}
Suppose that the Kronecker product $f\otimes g$  has an eigenvalue $\alpha\in K^\alg$ with algebraic multiplicity   
$\dim \sum_{l\geq 0}\Kernel(\alpha -f\otimes g)^l >mn-\min(m,n)$.
Then both factors $f$ and $g$ are homotheties.
\end{lemma}

\proof
It suffices to treat the case that $K$ is algebraically closed.    Set $I=\{1,\ldots, m\}$ and $J=\{1,\ldots,n\}$.
Consider the subset $S\subset I\times J$ comprising the pairs $(i,j)$ such that $\lambda_i\mu_j\neq \alpha$,
and let $E\subset V\otimes W$ be the generalized eigenspace with respect to  the eigenvalue $\alpha\in K$. 
Then $|S| = mn -\dim(E) < mn- (mn-\min(m,n))=\min(m,n)$.
In turn,  none of the projections $\pr_1:S\ra I$ and $\pr_2:S\ra J$ is surjective. In particular, there is
$1\leq r\leq m$ with $\alpha=\lambda_r\mu_j$ for all $j$. Using $\alpha\neq 0$ we obtain  $\mu_1=\ldots=\mu_n$,
hence $g$ is a homothety. By symmetry, the same holds for $f$.
\qed

\medskip
We also have used a general fact from non-abelian cohomology: Let $(X,e)$ be a pointed proper scheme over a ground field $F$,
endowed with a finite group $\Gamma\subset\Aut(X,e)$, and a $\Gamma$-stable closed subscheme $X_0$ containing $e$.
Then each $\Gamma$-torsor $P$ defines a twisted form $(X',X'_0,e')$ of the triple $(X,X_0,e)$, via
$$
X'=P\wedge^\Gamma X\quadand X'_0=P\wedge^\Gamma X_0,
$$ 
compare for example the discussion in \cite{Schroeer; Tziolas 2023}, Section 3.

\begin{lemma}
\mylabel{fact on twisting}
In the above situation, suppose that $\Aut_{(X_0,e)/F}$ is constant, and that the homomorphism $\Gamma\ra\Aut(X_0,e)$
is injective, with  normal image.
Let $P$ be a  $\Gamma$-torsor such that the resulting twisted form $(X',X'_0,e')$ satisfies $(X'_0,e')\simeq (X_0,e)$.
Then $(X',e)\simeq (X,e)$.
\end{lemma}

\proof
Write $G=\Aut(X_0,e)$ and   $\bar{G}=G/\Gamma$. The corresponding constant group schemes
sit in a short exact sequence $1\ra \Gamma_F\ra G_F\ra \bar{G}_F\ra 1$.
This is enough to yield in non-abelian cohomology a sequence of pointed sets
$$
G\lra \bar{G}\stackrel{\partial}{\lra} H^1(F,\Gamma_F)\lra H^1(F,G_F).
$$
By assumption, the class of $P$ in $H^1(F,\Gamma_F)$ maps to the class of the trivial torsor $G_F$  in $H^1(F,G_F)$.
As discussed in \cite{Hilario; Schroeer 2023}, Section 10  it lies in  the image of the coboundary map $\partial$.
The latter sends a point $a\in \bar{G}$ to its schematic fiber in $G_F$. All such fibers have a rational point,
thus $P$ is trivial, and hence $(X',e)\simeq (X,e)$.
\qed

\section{K3 surfaces over quadratic number fields}
\mylabel{K3 over quadratic}

Let $\QQ\subset F$ be a number field, $\O_F$ be the ring of integers, and   $S=\Spec(\O_F)$ the resulting Dedekind scheme. Write 
$$
\sigma_1,\ldots,\sigma_r\in S\quadand F_i^\wedge=\Frac(\O_{F,\sigma_i}^\wedge)\quadand F_i^\sh=\Frac(\O_{F,\sigma_i}^\sh)
$$
for the points of characteristic two and the ensuing   valued fields, the $F_i^\wedge$ being complete,
whereas the $F_i^\sh$  are strictly henselian.
We now examine the existence of elliptic curves  $E_F$ that are admissible for the Kummer construction,
and the ensuing question whether the resulting family of normal K3 surfaces $V=(E\times E)/G$ admits a simultaneous resolution
of singularities. As usual, $E\ra S$ is the Weierstra\ss{} model, and $G\subset\Aut_{E/S}$ denotes the effective model for the sign involution.
By Corollary \ref{potentially good reduction}, the admissible  $E_F$ have good reduction outside $\sigma_1,\ldots,\sigma_r\in S$,   potentially good reduction everywhere,
and the $j$-invariant is non-zero at each of the $\sigma_i$.
Let us summarize our general results of the preceding two sections for the situation at hand:

\begin{theorem}
\mylabel{sufficient conditions over number rings}
Notation as above. Suppose the elliptic curve $E_F$ is admissible for the Kummer construction.
Assume furthermore that for  each $1\leq i\leq r$ the following holds:
\begin{enumerate}
\item 
If the elliptic curve $E_F$ has good reduction at $\sigma_i\in S$,   the group scheme $E_F[2]$ becomes constant over $F_i^\wedge$.
\item 
If $E_F$ has bad reduction at $\sigma_i$, it acquires good reduction over some quadratic extension of $F_i^\sh$,
the effective model $G_{\sigma_i}$ is multiplicative, and
the field $F$ contains a primitive third root of unity.
\end{enumerate}
Then the family of normal K3 surfaces $V=(E\times E)/G$ admits  simultaneous resolutions of singularities,
and in particular $\shM_\Kthree(\O_F)\neq\varnothing$.
\end{theorem}

\proof
This follows from Theorem \ref{simultaneous partial resolution} and Theorem \ref{sufficient condition for resolution}.
\qed

\medskip
Parsing the  literature on elliptic curves over number fields, one indeed finds examples. Here we consider
quadratic number fields, and will turn to $S_3$-number fields in the next section.
 
Recall that every quadratic extension $\QQ\subset L$ is cyclic, and takes the form $L=\QQ(\sqrt{m})$
for some unique square-free integer $m\neq 0,1$. Up to isomorphism, the field $L$ is determined by its discriminant
$$
d_L=\begin{cases}
m	& \text{if $m\equiv 1$ modulo 4;}\\
4m	& \text{else.}
\end{cases}
$$
The occurring values, sometimes called ``fundamental discriminants'',   are  precisely the numbers
$d_L=\epsilon \left(\frac{-1}{p_1\ldots p_r}\right)\cdot 2^\nu p_1\ldots p_r $ 
for  pairwise different odd primes $p_i$ and $\nu\in\{0,2,3\}$ and $\epsilon=\pm 1$,
subject to the   constraint   $\nu=2\Rightarrow \epsilon=-1$. Note that $d_L=12$ is impossible.

Consider the following eight elliptic curves $E_F:y^2=x^3+a_2x^2+a_4x$
over certain quadratic number fields $F=\QQ(\sqrt{m})$:
$$
\begin{array}{llllll}
\toprule
d_F	& m	&  \epsilon	& a_2			& a_4			& j\\
\toprule
28	& 7	&  8+3\sqrt{7}	& -(1+2\epsilon^2)	& 16\epsilon^3  	& 255^3\\
	&   	&  		& -(1+2\epsilon^{-2})	& 16\epsilon^{-3}	& \\
\midrule
41	& 41	& 32+5\sqrt{41}	& (3\epsilon-1)/2	& (\epsilon^2-\epsilon)/2	& (\epsilon-16)^3/\epsilon\\
	& 	& 		& (-3\epsilon^{-1}-1)/2	& (\epsilon^{-2}+\epsilon^{-1})/2\\
\midrule
65	& 65	& 8+\sqrt{65}	& 2\epsilon^2-1		& 16\epsilon^3		& 257^3\\
	& 	& 		& 10\epsilon^2-5	& 400\epsilon^3	\\
	& 	& 		& 8\epsilon+1		& 16\epsilon^2		& 17^3\\
	& 	& 		& 40\epsilon+5		& 400\epsilon^2\\
\bottomrule
\end{array}
$$
Here $d_F$ is the discriminant and $\epsilon\in \O_F^\times$ is a fundamental unit.
According to \cite{Comalada 1990}, Theorem 2, these $E_F$ have good reduction everywhere
and constant group scheme $E_F[2]$. Conversely,  every   elliptic curve over a quadratic number field
with these properties belongs to the list. Clearly,   $j(\sigma)\neq 0$ for each point $\sigma\in S$ of characteristic two.

Let $F$ and a pair $E_F$, $E'_F$ be as in the table, with resulting  Weierstra\ss{} models $E$ and $E'$.
According to Theorem \ref{sufficient conditions over number rings}, the family of normal K3 surfaces $V=(E\times E')/G$ admits a simultaneous resolution
of singularities $X\ra V$. This shows:

\begin{theorem}
\mylabel{k3 over   quadratic fields}
We have $\shM_\Kthree(\O_F)\neq\varnothing$ for the quadratic number fields $F$ with  discriminant $d_F\in\{28,41,65\}$.
\end{theorem}

We actually get more:
Write $\shM_\Enr\ra (\Aff/\ZZ)$ for the \emph{stack of Enriques surfaces}. The objects are triples $(R,Y,f)$ where $R$
is a ring, $Y$ is an algebraic space, and $f:Y\ra\Spec(R)$ is a proper flat morphism of  finite presentation
whose geometric fibers are Enriques surfaces.

\begin{theorem}
\mylabel{enriques over quadratic fields}
We have $\shM_\Enr(\O_F)\neq\varnothing$ for the quadratic number fields $F$  with  discriminant $d_F\in\{28,41,65\}$.
\end{theorem}

\proof
Let $\zeta_1,\ldots,\zeta_3\in E(F)$ be the three points of order two.
At each point $\sigma\in S$ with residue field $k=\kappa(\sigma)$ of characteristic two,
the kernel for the specialization map $E(F)\ra E(k)$
contains exactly one of the $\zeta_i$. Since there are at most two such $\sigma$,
we find some $\zeta=\zeta_j$ such that $\zeta\in E(\O_F)$ has order two in all fibers.
Likewise, we choose such  $\zeta'\in E'(\O_F)$. 
On the family of abelian surfaces $A=E\times E'$, this defines an involution   $(x,x')\mapsto (x+\zeta,-x'+\zeta')$.
These are often called \emph{Lieberman involutions}.
The resulting action of $H=(\ZZ/2\ZZ)_S$ is free, because it is free on the first factor, and commutes
with the sign involution, since $\zeta$ and $\zeta'$ coincide with their negatives. 
In turn, we get an induced $H$-action on the family $Y=A/G$
of normal K3 surfaces.

Let $Y=\Bl_Z(V)$ be the blowing-up as in \eqref{first blowing-up curve family}, which by Proposition \ref{center with embedded components}
is fiberwise given by blowing-up the singularities. Our $H$-action on $V$ uniquely extends to $Y$,
because the center $Z\subset V$ is $H$-stable.
In turn, the exceptional divisors  $\Delta_1,\ldots,\Delta_{16}$ on the generic fiber $Y_\eta$ are $H_s$-stable,
and so are their schematic images in $Y$.
By Proposition \ref{sixteen blowing-ups}, their successive blowing-ups    yield a simultaneous resolution of singularities $X\ra Y$,
and the $H$-action again extends uniquely to $X$.

We claim that the $H$-action on $X$ is free. 
The induced action on $\Exc(X/V)$ is free, because this holds for  $\Sing(V/S)$.
It remains to verify that the action on $\Reg(V/S)$ is free.
This may be checked fiberwise, with geometric points.
Suppose   that  $(x,x')$ is such a  point on $A_s$ that maps to a $H_s$-fixed point of  $\Reg(V/S)_s$.
Then $-(x,x')=(x+\zeta,-x'+\zeta')$, and thus $\zeta'=0$, contradiction. In turn, the free quotient $X/H$ yields
the desired family of Enriques surfaces.
\qed

\medskip
We close this section with the following observation:

\begin{proposition}
\mylabel{no curve over small fields}
There are no elliptic curves over $\QQ$ or $F=\QQ(\sqrt{m})$ for $m\in\{-2,-1\}$ that are admissible for the Kummer construction.
\end{proposition}

\proof 
Ogg  classified the elliptic curves $E_\QQ$   whose discriminant take the form $\Delta=\pm 2^v$,
there are twenty-four isomorphism classes (\cite{Ogg 1966}, Table 1). In each case $E$ has a Weierstra\ss{} equation of the form
$y^2=x^3+a_2x^2+a_4x$ with  $a_2\in 2\ZZ$ and $\nu\geq 1$. Thus  both $\Delta$ and  $a_2a_4+a_6$ belong to $2\ZZ$,
and Proposition \ref{no intersection} tells us that $E_\QQ$ is not admissible.

We next examine  $F=\QQ(\sqrt{-1})$.  Write $i=\sqrt{-1}$ for the imaginary number.
The ring of integers is given by $\O_F=\ZZ[i]$, with group of units $\O_F^\times=\{i^n\mid 0\leq n\leq 3\}$, and  uniformizer 
$\pi=1+i$ at  the point $\sigma\in S$ of characteristic two. 
Pinch  classified the $E_F$ whose discriminant takes the form $i^n\pi^\nu$,
now  there are sixty-four isomorphism classes (\cite{Pinch 1984}, Table 2). In all cases $\nu\geq 1 0$, so $\Delta\in\pi\O_F$.
All but one have a   Weierstra\ss{} equation of the form $y^2=x^3+a_2x^2+a_4x+a_6$.
Going through these cases and using that both $2$ and $1+i$ belong to the maximal ideal $\pi\O_F$, one finds $a_2a_4+a_6\in \pi\O_F$.
These $E_F$ are not admissible by Proposition \ref{no intersection}.
The remaining case has
$$
y^2+(1+i)xy=x^3+ix^2+2x+3i.
$$
Using \eqref{relation b-values a-values} we compute  $b_2a_4+b_6 =   24i$, thus  $\val(b_2a_4+b_6)=6$. 
In light of  $2=i(1+i)^2$ one gets $\val(2,a_1,a_3)=1$,
so   $E_F$ is not admissible for the Kummer construction by Proposition \ref{no intersection}.

Finally consider the imaginary quadratic number field $F=\QQ(\sqrt{-2})$.
Now the ring of integers is $\O_F=\ZZ[\sqrt{-2}]$, with group of units $\O_F^\times=\{\pm 1\}$,
and uniformizer $\pi=\sqrt{-2}$ at the point $\sigma\in S$ of characteristic two.
Pinch also classified the $E_F$ whose discriminant takes the form $\pm\pi^\nu$,
now there are forty cases. All but two have Weierstra\ss{} equation of the form 
$y^2=x^3+a_2x^2+a_4x+a_6$.
Going through these cases,  one finds $a_2a_4+a_6\in \pi\O_F$, so these $E_F$ are not admissible by Proposition \ref{no intersection}.
The remaining two cases have Weierstra\ss{} equations
$$
y^2+\pi xy=x^3-x^2-2x+3\quadand y^2+\pi xy+\pi y=x^3-x^2-x.
$$
Using \eqref{relation b-values a-values} we compute the respective values $b_2a_4+b_6= 24 $ and $b_2a_4+b_6= 4$.
In both cases we see  $\val(b_2a_4+b_6)\geq 4$ and $\val(2,a_1,a_3)=1$.
Again by Proposition  \ref{no intersection}, the $E_F$ are not admissible for the Kummer construction.
\qed

\section{K3 surfaces over  \texorpdfstring{$S_3$}{S3}-number fields}
\mylabel{K3 over s3}

Recall that an \emph{$S_3$-number field} is a finite Galois extension $\QQ\subset F$ whose Galois group is isomorphic to the
symmetric group $S_3$. The goal of this final section is to establish:

\begin{theorem}
\mylabel{k3 over s3 fields}
Let $\QQ\subset F$ be the $S_3$-number field arising as Galois closure from 
a cubic number field $K$  with discriminant $d_K=-3f^2$ for some even $f\geq 1$.
Then $\shM_\Kthree(\O_F)\neq \varnothing$.
\end{theorem}

The proof appears at the end of this section, after we    have introduced  
certain elliptic curves and  reviewed the relevant facts about $S_3$-number fields.

Consider the  imaginary quadratic number field  $L=\QQ(\sqrt{-3})=\QQ(e^{2\pi i/3})$ and set $\omega=e^{2\pi i/6}=(1+\sqrt{-3})/2$.
Note that the ring of integers is $\O_L=\ZZ[\omega]$, the sixth root of unity $\omega\in  \O_L^\times$ is a generator,
the prime  $p=2$ is inert,  and the resulting $\primid=2\O_L$ has residue field $\FF_4$.   Consider the two elliptic curves 
\begin{equation}
\label{pinch equations}
E_L:\quad y^2=x^3 \pm(1+\omega)x^2 +\omega x.
\end{equation} 
These are the first entries in \cite{Pinch 1984}, Table 4, where Pinch   classified all elliptic curves over $L$
that have good reduction outside $p=2$. 

\begin{proposition}
\mylabel{properties pinch curves}
Up to isomorphism,  \eqref{pinch equations} are the only elliptic curves over $L$
that are admissible for the Kummer construction. For both of them, 
the effective model $G$ of the sign involution is isomorphic to $\mu_{2,S}$. Moreover,
the invariant is $j=0$, the discriminant is $\Delta=2^4$,       the reduction has Kodaira symbol $\II$,
and the finite \'etale group scheme $E_L[2]$ becomes constant over the henselization $\O_L^h$ with respect to $\primid=2\O_L$.
\end{proposition}

\proof
By the classification of Pinch in loc.\ cit., there are fifty-four cases of elliptic curves over $L$
with discriminant $\Delta=\pm\omega^{\pm 1}2^\nu$, each given by a Weierstra\ss{} equation of the form $y^2=x^3+a_2x^2+a_4x + a_6$.
We thus have $2_d=1$; so by Proposition \ref{no intersection} admissibility is equivalent to $a_2a_4+a_6\not\equiv 0$ modulo $2\O_L$.

For \eqref{pinch equations} we have $a_2a_4+a_6=\pm(1+\omega)\omega\equiv 1$ modulo $2\O_L$, so our
$E_L$ are admissible for the Kummer construction. One directly computes the discriminant and $j$-invariant,
and the Tate Algorithm \cite{Tate 1975} reveals that the Kodaira symbol is $\II$.
Since $2\in\O_L$ is prime, the Tate--Oort classification (\cite{Tate; Oort 1970}, Theorem 2) ensures
that $G$ is isomorphic to either $(\ZZ/2\ZZ)_S$ or $\mu_{2,S}$.
By Proposition \ref{fixed scheme etale} the effective model $G$ of the sign involution is multiplicative,
and we get $G=\mu_{2,S}$. The group scheme $E_L[2]$ is given by the equations $x=0$ and $x^2\pm (1+\omega)x+\omega=0$.
Over the residue field $\FF_4$, the latter factors into $(x-\omega)(x-1)$, and the assertion
follows from Hensel's Lemma.

For all remaining Weierstra\ss{} equations from Pinch's  table, one checks  $a_2a_4+a_6\equiv 0$ modulo $2\O_L$,
so there are no  further admissible curves.   
\qed

\medskip
Choose a separable closure for the residue field $\O_L/\primid=\FF_4$, form the ensuing
strictly henselian ring $\O_L^\sh$ with respect to $\primid=2\O_L$, and let $L^\sh$ be its field of fractions.
According to \cite{Serre; Tate 1968}, Corollary 3 for Theorem 2,
 there is a smallest extension $L^\sh\subset L^\good$ over which $E_{L^\sh}$ acquires
good reduction, and this extension is Galois.
Moreover, the  Galois group is isomorphic to a  subgroup of $\Sp_2(\FF_3)=\SL_2(\FF_3)$, a group of order twenty-four.
It can also be seen as the semidirect product $Q\rtimes C_3$ of the quaternion group $Q=\{\pm E,\pm I,\pm J,\pm K\}$
by a cyclic group of order three.

\begin{proposition}
\mylabel{galois group for good reduction}
The   group $\Gal(L^\good/L^\sh)$ is  isomorphic to $C_2\times C_3$.
\end{proposition}

\proof
We have $j=0$, hence $c_4=0$, and furthermore $\Delta=2^4$. 
The   ramification index for $p=2$ in $\O_L$ is $e=1$. We now apply \cite{Kraus 1990}, Theorem 2:
Since  $3\val(c_4)=\infty\geq 12e+\val(\Delta)$ and $3\nmid\val(\Delta)$ and the Kodaira symbol is not $\IV$ or $\IV^*$,
the Galois group has order six. It must be $C_2\times C_3$, because $-E\in Q$ is the only element of order two
in $Q\rtimes C_3$.
\qed

\medskip
We see that the cyclic extension $L^\sh\subset L^\good$ contains four intermediate fields.
Of particular interest is what we call   $L^\bad$, the intermediate field with $[L^\bad:L^\sh]=3$.

\begin{proposition}
\mylabel{reduction type after first base change}
For the base-change $E\otimes {L^\bad}$, the Weierstra\ss{} equations \eqref{pinch equations} remain minimal,
but the    Kodaira symbol changes to  $\II^*$. Moreover, it acquires good reduction over the quadratic extension
$L^\bad\subset L^\good$.
\end{proposition}

\proof
Write  $R^\bad \subset L^\bad $ for the integral closure of $R^\sh=\O_L^\sh$.
First note that by  construction, $E\otimes L^\bad$ has bad reduction, that $R^\sh\subset R^\bad$ is totally ramified,
and that the  Weierstra\ss{} equation \eqref{pinch equations} has $\val_{R^\sh}(\Delta)=4$ and $\val_{R^\bad}(\Delta)=12$.
It follows that the equations remain  minimal. The extension $L^\bad\subset L^\good$ has degree two.
By definition, good reduction occurs over this extension.

Recall that for additive reduction, Ogg's Formula (\cite{Ogg 1967}, see also the discussion in \cite{Schroeer 2021a}) takes the form 
$\val(\Delta)= 2 + \delta+ (m-1)$, where $m\geq 1$ is the number of irreducible components in the closed fiber of the minimal regular model,
and $\delta$ is the wild part for the conductor of the Galois representation on   $\ell$-torsion points in the generic fiber.
Over $R^\sh$, this becomes $4=2+\delta+(1-1)$, hence $\delta=2$.
The extension $R^\sh\subset R^\bad$ is tamely ramified, so the $\delta$ does not change when passing to $R^\bad$.
Over this   ring, Ogg's Formula yields $12=2+2+(m'-1)$, so $m'=9$. The only possible Kodaira symbols  are 
$\II^*$ and $\I_4^*$. The latter is impossible with $c_4=j=0$, according to \cite{Papadopoulos 1993}, Table V. 
\qed

\medskip
Recall that over a given field $k$,  the Galois group of an irreducible separable polynomial $P(T)$ of degree three
is isomorphic to $S_3$ if and only if $\disc(P)\in k^\times $ is not a square. By the Galois Correspondence,
the  resulting splitting field $k\subset F$
then contains three   cubic fields $K_1,K_2,K_3$ and another quadratic field  $L$. Moreover, $F$ coincides
with  the Galois closure for each of the $K=K_i$.
For $k=\QQ$  the finite extension $\QQ\subset F$ is called an \emph{$S_3$-number field}. 
 
Suppose $k\subset K$ is a separable cubic extension.
Such extensions are twisted forms of the product ring $k\times k\times k$,
hence the group scheme $\Aut_{K/k}$ is an inner form of the constant group scheme $(S_3)_{k}$, as described in 
\cite{Schroeer; Tziolas 2023}, Lemma 3.1.  Since $A_n\subset S_n$, $n\geq 3$ is the commutator subgroup,
we get an induced twisted form $\Aut'_{K/k}$ of  $(A_3)_k$. The underlying scheme  takes the form 
$\Aut'_{K/k}=\{e\}\cup \Spec(L)$ for some \'etale $k$-algebra $L$ of degree two.
 
\begin{lemma}
\mylabel{non-normal cubic extensions}
The separable cubic  extension $k\subset K$ is non-normal if and only if  the \'etale algebra 
$L$ is a field. In this case, the Galois closure is $K\otimes_kL$, and the following are equivalent:
\begin{enumerate}
\item
The cubic extension takes the form  $K=k(\sqrt[3]{\alpha})$ for some  $\alpha\in k^\times$.
\item 
The   extension $k\subset L$ is obtained by adjoining a primitive third root of unity.
\item 
The group scheme $\Aut'_{K/k}$ is isomorphic to $\mu_3$.
\end{enumerate}
In characteristic $p\neq 2$, we actually have $L=k(\sqrt{-3})$.
\end{lemma}

\proof
If not a field, $L$ must be isomorphic to $k\times k$. It follows that $\Aut(K/k)$ contains a subgroup of order $[K:k]=3$,
hence the   extension $k\subset K$ is Galois. Conversely, if the separable extension $k\subset K$ is normal,
the Galois group must be cyclic of order three. But $\Aut'_{K/k}$ is the only subgroup scheme
of order three inside  $\Aut_{K/k}$, and it follows that $L=k\times k$. This establishes the first assertion.

Suppose now that $K$ is non-normal, so $L$ is  a field. For degree reasons, $\tK=K\otimes_kL$ remains a field,
and by the preceding paragraph the cubic extension  $L\subset\tK$ is cyclic. Moreover, $K\subset\tK$ is separable
of degree two. It follows that $\Aut(\tK/F)$ contains elements of order three and two, hence $k\subset\tK$
is a Galois extension of degree six. The Galois group is non-abelian, because the intermediate extension $k\subset K$ is non-normal,
and it follows that $\tK$ is the Galois closure.

It remains to establish the equivalence of (i)--(iii). 
First note that each of the three conditions implies $p\neq 3$, and we thus may disallow this characteristic.
Suppose first  $K=k(\sqrt[3]{\alpha})$. Then  the Galois closure $\tK$ is obtained by adjoining a primitive third root of unity $\zeta$,
and thus $L=K(\zeta)$.
Then $\omega=-\zeta$ has minimal polynomial $Q(T)=T^2-T+1$, and consequently $(2\omega+1)^2=-3$.
\qed

\medskip
We now restrict to   $k=\QQ$, that is, \emph{cubic number fields} $\QQ\subset K$. For more details we   refer to the seminal work 
of Hasse \cite{Hasse 1930} and  the monographs of Delone and Faddeev \cite{Delone; Faddeev 1964}, Cohen \cite{Cohen 1993},
and Hambleton and Williams \cite{Hambleton; Williams 2018}. 
The  \emph{pure cubic fields}, which can be written as  $K=\QQ(\sqrt[3]{m})$, are of particular interest.

\begin{lemma}
\mylabel{pure cubic number fields}
A cubic  number field $\QQ\subset K$ is pure if and only if the discriminant takes the form $d_K=-3f^2$. 
Then $K$ is  non-normal, with quadratic field $L=\QQ(\sqrt{-3})$. 
\end{lemma}

\proof 
The first assertion is well-known, an  argument appears  in \cite{Guardia; Pedret 2025}, Section 7 
or \cite{Bhargava; Shnidman 2014}, Lemma 33.
For the second assertion, suppose the equivalent conditions hold. Then $K$ has  a real embedding,
and its Galois closure contains  a primitive third root of unity.  Hence it   is non-normal, with $L=\QQ(\sqrt{-3})=\QQ(e^{2\pi i/3})$.
\qed

\medskip
\emph{Proof of Theorem \ref{k3 over s3  fields}.} We start with the quadratic number field
$L=\QQ(\sqrt{-3})=\QQ(e^{2\pi i/3})$, set $\omega=e^{2\pi i/6}=(1+\sqrt{-3})/2$, 
and  fix one of the elliptic curves
$$
E_L:\, y^2=x^3 \pm(\omega+1)x^2+\omega x. 
$$
We saw in Proposition \ref{properties pinch curves} that $E_L$ is  admissible for the Kummer construction,
with $G=\mu_2$ as    effective model for the sign involution.
In turn, the categorical quotient $V=(E\times E)/G$ yields  a family of normal K3 surfaces over $S=\Spec(\O_L)$.
According to Lemma \ref{center with embedded components}, there is a simultaneous partial resolution $Y\ra V$
for all singularities stemming from the fixed points.
Now $\Sing(Y/S)=\{y_\crit\}$ is a singleton,  which is a $D_4$-singularity in its  fiber over the residue field
$\O_L/2\O_L=\FF_4$.

The given cubic number field $\QQ\subset K$ has discriminant $d_K=-3f^2$ with some even $f\geq 1$.
By Lemma \ref{pure cubic number fields} we have $K=\QQ(\sqrt[3]{m})$ for some $m\geq 1$, 
and the Galois closure is given by $F=K\otimes L$. The extension $L\subset F$ is cyclic of degree three.
The assumption on $f$ ensures that $\O_L\subset\O_F$ is totally ramified over $\primid=2\O_L$.

Fix a separable closure for $\O_L/\primid=\FF_4$, with ensuing strict henselization $\O_L^\sh$
and field of fractions $L^\sh=\Frac(\O_L^\sh)$.
Also choose $F^\sep$ and an embedding $L^\sep\subset F^\sep$.
We now want to apply Theorem \ref{sufficient condition for resolution}, and have to verify its  assumptions (i)--(iii).
Obviously, our field $L=\QQ(e^{2\pi i/6})$ contains a primitive third root of unity,
and the purely inseparable extension $\FF_4\subset\kappa(v_\crit)$ is an equality.
Over $p=2$, the extension $\ZZ\subset\O_L$ is \'etale and $\ZZ\subset\O_K$ is totally ramified, hence $\O_L\subset\O_F$
is totally ramified over $\primid=2\O_L$.
In turn, the induced extension $L^\sh\otimes_LF$ coincides with $L^\sh\subset L^\good$.
By Proposition \ref{reduction type after first base change}, 
the elliptic curve $E\otimes F^\sh$ acquires good reduction over some quadratic extensions.
So the theorem applies, and the base-change $V\otimes\O_{\tK}$ admits a simultaneous minimal resolution of singularities.
\qed



\end{document}